\documentclass[a4paper,10pt]{article}
\usepackage{amssymb,amsmath,amsthm,amsfonts}
\usepackage{bookmark}
\usepackage{mathrsfs}
\usepackage{multirow}
\usepackage{graphicx}
\usepackage{authblk}
\usepackage{indentfirst}
\usepackage{multicol}
\usepackage{tabu}
\usepackage{url}
\usepackage{fancyhdr}
\usepackage[numbers]{natbib}
\usepackage[all]{xy}
\usepackage{mathtools}
\usepackage[english]{babel}
\usepackage{colortbl}
\usepackage{caption}
\usepackage{hyperref}\usepackage{subfigure}
\usepackage{tikz}\usepackage{float}\usepackage{CJK}
\usepackage{xcolor}\usepackage{mdframed}\usepackage{footmisc}

\newtheorem{example}{Example}[section]

\newcommand{\im}{{\mathrm{Im}\hspace{0.1em}}}

\newcommand{\rank}{{\mathrm{rank}\hspace{0.1em}}}
\newcommand{\Inf}{{\mathrm{Inf}\hspace{0.1em}}}
\newcommand{\Sup}{{\mathrm{Sup}\hspace{0.1em}}}

\definecolor{blue}{rgb}{0.0, 0.6, 1.0}
\newcommand{\blue}[1]{{\color{blue} #1}}
\definecolor{red}{rgb}{1.0, 0.0, 0.0}
\newcommand{\red}[1]{{\color{red} #1}}
\newcommand{\gray}[1]{{\color{gray} #1}}
\newcommand{\orange}[1]{{\color{orange} #1}}
\definecolor{green}{rgb}{0.00,0.60,0}
\newcommand{\green}[1]{{\color{green} #1}}

\usepackage{xr}
\makeatletter
\newcommand*{\addFileDependency}[1]{
	\typeout{(#1)}
	\@addtofilelist{#1}
	\IfFileExists{#1}{}{\typeout{No file #1.}}
}
\makeatother

\title{ChatGPT for Computational Topology}

\author[1,2]{Jian Liu \thanks{Corresponding author: liujian2@mail.nankai.edu.cn}}
\author[2]{Li Shen}
\author[2,3,4]{Guo-Wei Wei \thanks{Corresponding author: weig@msu.edu}}
\affil[1]{Mathematical Science Research Center, Chongqing University of Technology, Chongqing 400054, China}
\affil[2]{Department of Mathematics, Michigan State University, MI, 48824, USA}
\affil[3]{Department of Electrical and Computer Engineering, Michigan State University, MI 48824, USA}
\affil[4]{Department of Biochemistry and Molecular Biology, Michigan State University, MI 48824, USA}

\makeatletter
\renewcommand*{\@fnsymbol}[1]{\ensuremath{\ifcase#1\or \dagger\or *\or *\or
		\mathsection\or \else\@ctrerr\fi}}
\makeatother
\date{}

\begin{document}
	\begin{CJK*}{GBK}{kai}
		\CJKtilde
		\maketitle
		
		\paragraph{Abstract}
		
	ChatGPT represents a significant milestone in the field of artificial intelligence (AI), finding widespread applications across diverse domains. However, its effectiveness in mathematical contexts has been somewhat constrained by its susceptibility to conceptual errors. Concurrently, topological data analysis (TDA), a relatively new discipline, has garnered substantial interest in recent years. Nonetheless, the advancement of TDA is impeded by the limited understanding of computational algorithms and coding proficiency among theoreticians. This work endeavors to bridge the gap between theoretical topological concepts and their practical implementation in computational topology through the utilization of ChatGPT. We showcase how a pure theoretician, devoid of computational experience and coding skills, can effectively transform mathematical formulations and concepts into functional codes for computational topology with the assistance of ChatGPT.
Our strategy outlines a productive process wherein a mathematician trains ChatGPT on pure mathematical concepts, steers ChatGPT towards generating computational topology codes, and subsequently validates  the generated codes using established examples.
Our specific case studies encompass the computation of Betti numbers, Laplacian matrices, and Dirac matrices for simplicial complexes, as well as the persistence of various homologies and Laplacians. Furthermore, we explore the application of ChatGPT in computing recently developed topological theories for hypergraphs and digraphs, as well as the persistent harmonic space, which has not been computed in the literature, to the best of our knowledge. This work serves as an initial step towards effectively transforming pure mathematical theories into practical computational tools, with the ultimate goal of enabling real applications across diverse fields.




	     \newpage
		 {\setcounter{tocdepth}{4} \tableofcontents}
	     \newpage
		
		\section{Introduction}\label{section:introduction}

		Mathematical theory guides the advancement of applied sciences, while the progress of applied sciences, in turn, provides direction and motivation for mathematical exploration. However, the development of these disciplines can inadvertently erect barriers between them. A theoretical mathematician often finds it challenging to divert their focus toward translating their theoretical research into practical programs that can be utilized by researchers in the applied sciences for numerical experiments. Conversely, an expert in applied sciences may struggle to grasp the intricacies of advanced mathematical theories and harness this knowledge to optimize their methods and models. This interplay between theory and application highlights the need for bridging the gap and fostering collaboration between these two seemingly disparate domains.

		Recently, ChatGPT has excelled in various domains and has been driving advancements in artificial intelligence (AI). For example, ChatGPT is being used to assist in drug discovery and development \cite{sharma2023chatgpt,wang2023chatbots}. Harnessing ChatGPT to assist theoretical mathematicians in developing algorithms and codes, even without proficiency in computer languages, can be a game-changer. This approach addresses the critical need for collaboration between theoretical and applied scientists, as it facilitates knowledge exchange and interdisciplinary research. One promising avenue is to use natural language interfaces to interact with ChatGPT. Theoretical mathematicians can engage in conversations with ChatGPT, train ChatGPT with correct mathematical concepts or theories, and describe their computational objectives in plain language. ChatGPT can then assist them in translating these theories into algorithms and codes, without requiring an in-depth understanding of specific programming languages and algorithms by theoreticians. It can provide code snippets, algorithm frameworks, and executable scripts to help theoretical mathematicians initiate practical computational work.
		
	Topology is one of the most profound and challenging areas in modern mathematics. In recent years, topological data analysis (TDA) has gained significant attention in applied fields such as computer science, molecular biology and biophysics, materials science, medical imaging, etc. An increasing number of scholars and researchers are dedicating themselves to develop the theories on TDA \cite{carlsson2009topology,carlsson2004persistence,edelsbrunner2002topological,zomorodian2004computing}. The computational topology is an emerging field in mathematics \cite{edelsbrunner2022computational,kaczynski2004computational}. However, the lack of algorithm understanding and coding skills limits mathematicians' success in computational topology.
In this work, we aim to leverage ChatGPT to close the gap between pure mathematical theories and practical codes for computational topology.

		Some traditional perspectives  argue that ChatGPT may not be suitable for tasks involving complex mathematical concepts and intricate logical reasoning due to its limited depth in mathematics. Many limitations stem from several factors. First, ChatGPT's understanding of mathematics is based on pre-existing knowledge up to its last training data, which means it may not be aware of recent developments or specialized mathematical domains. Second, while it can perform various mathematical calculations and explain some fundamental concepts, it lacks the depth of understanding and rigor that specialized mathematicians or dedicated mathematical software can provide. Additionally, ChatGPT may struggle with highly specialized or research-level problems that require nuanced domain knowledge and precise mathematical proofs. Consequently, when employing ChatGPT for coding tasks, users must supplement its capabilities with their own expertise in mathematics and employ specific strategies to supervise and guide ChatGPT in generating accurate codes.
		
		When using ChatGPT, it is essential to keep the following considerations in mind:
		\begin{itemize}
			\item \textbf{Answers May Not Be Definitive:} ChatGPT's responses to questions are not guaranteed to be correct or entirely satisfactory. Its answers are based on patterns and information in its training data, and it may not always provide accurate or complete responses.
			\item \textbf{Variability in Responses:} The same question posed to ChatGPT in different instances or chat interfaces may yield different answers. This variability can occur due to the model's inherent randomness and sensitivity to the input phrasing.
			\item \textbf{Responses Evolve with Interaction:} In a single chat session, ChatGPT's answers can evolve and become closer to what the user desires as the conversation progresses. It adapts to the user's guidance and context.
		\end{itemize}
During the interaction with ChatGPT, we always assume that the questioners are familiar the related mathematical theories but lacks coding experience, skill, and understanding.

		In the next section, we will recall the definitions and notations for simplicial complexes and Betti numbers. We   demonstrate how to obtain practical codes for computing the Betti numbers, Laplacian, and Dirac   of a simplicial complex using ChatGPT.  We intend to use the Python programming language for this purpose, as Python is widely adopted in the scientific and data analysis communities. Also, Python offers a rich ecosystem of libraries and tools for scientific computing and data analysis which can significantly simplify the implementation of mathematical algorithms.
		In Section \ref{section:hypergraph_laplacian}, we   show how to use ChatGPT to discover the codes for computing the Betti numbers and Laplacian matrices for hypergraphs and digraphs.
		Section \ref{section:persistence} provides a note on the persistence of homology, Laplacians, and harmonic space.
In Section \ref{section:application}, we present some application  examples. Finally, we summarize the main works of this paper and outline further perspectives about ChatGPT.
ChapGPT version 3.5 was used in this work.
		
\section{ChatGPT for homology, Laplacian, and Dirac}\label{section:homology}

		Persistent homology is one of the fundamental theoretical tools in TDA, which can describe the multi-scale information of data sets \cite{carlsson2004persistence,edelsbrunner2002topological,edelsbrunner2008persistent,zomorodian2004computing}. Persistent Laplacians, including persistent spectral graph \cite{wang2020persistent} and evolutionary de Rham-Hodge theory \cite{chen2019evolutionary}, were introduced in 2019. Persistent Laplacians not only provide geometric information of data sets in terms of topological invariants as  persistent homology does, but also capture homotopic shape evolution that cannot be detected by persistent homology \cite{wang2020persistent,liu2023algebraic,memoli2022persistent}. Their applications in real-world problems show remarkable performance  \cite{chen2022persistent,meng2021persistent}. Furthermore, the persistent Dirac has also been introduced as a new computational tool in TDA \cite{ameneyro2022quantum,wee2023persistent}.
In this section, we   review the concepts of simplicial complexes, homology, and Betti numbers. Then we   showcase the use of  ChatGPT for computing Betti numbers from a mathematical perspective. Moreover, we  instruct ChatGPT to produce the codes for computing Dirac and Laplacian operators.

Our questions to ChatGPT and the design of corresponding codes are based on our mathematical background on homology theory. In this work, the ground field is always assumed to be $\mathbb{K}$. In particular, when we consider the Laplacian, the ground field is taken to be the real number field $\mathbb{K}=\mathbb{R}$.
		
		\subsection{Simplicial complexes and Betti numbers}
		
Computers are more adept at handling and processing discrete objects. As a result, researchers often focus their efforts on investigating simplicial complexes or their discrete representations, known as abstract simplicial complexes.
		
		Let $\mathcal{K}$ be an abstract simplicial complex on a finite nonempty ordered set $V$. Let $C_{p}(\mathcal{K})$ be the vector space with basis given by the $p$-simplices of $\mathcal{K}$.
		We have the chain complex $C_{\ast}(\mathcal{K})=(C_{p}(\mathcal{K}))_{p\geq 0}$. The $p$-th boundary operator $d_{p}:C_{p}(\mathcal{K})\to C_{p-1}(\mathcal{K})$ is given by
		\begin{equation*}
			d_{p}[v_{0},\dots,v_{p}]=\sum\limits_{i=0}^{p}(-1)^{i}[v_{0},\dots,\hat{v_{i}},\dots,v_{p}].
		\end{equation*}
		Here, $[v_{0},\dots,v_{p}]$ is a $p$-simplex of $\mathcal{K}$, and $\hat{v_{i}}$ means omission of $v_{i}$. Choose the basis of $C_{p}(\mathcal{K})$ as the lexicographic order of the $p$-simplices. Then we can obtain the representation matrix $B_{p}$ of the boundary operator $d_{p}$. Then the $p$-th Betti number of $\mathcal{K}$ is given by
		\begin{equation}\label{equ:betti}
			\beta_{p}=n_{p}-\rank B_{p}-\rank B_{p+1}.
		\end{equation}
		The Betti numbers for a simplicial complex can be computed by the above formula.
		
\begin{figure}[H]
			\centering
			\begin{mdframed}[backgroundcolor=lightgray, linewidth=0pt]
				\noindent Q: Given an abstract simplicial complex, can you give me a Python code for computing Betti numbers for the simplicial complex?
			\end{mdframed}
			\begin{mdframed}[backgroundcolor=orange!30, linewidth=0pt]
				\noindent \includegraphics[width=0.15in]{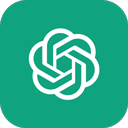} 
				Certainly! You can compute Betti numbers for an abstract simplicial complex using various Python libraries, such as \textbf{`networkx'} and \textbf{`scipy'}. First, make sure you have these libraries installed:
               \begin{mdframed}[backgroundcolor=black, linewidth=0pt]
					\textcolor[rgb]{1.00,1.00,1.00}	{
						\texttt{pip install networkx scipy
					}}
				\end{mdframed}
                Here's a Python code example to compute Betti numbers for a simplicial complex:
				\begin{mdframed}[backgroundcolor=black, linewidth=0pt]
					\textcolor[rgb]{1.00,1.00,1.00}{ \texttt{\blue{import} networkx \blue{as} nx\\
							\blue{import} scipy.linalg\\
							\blue{def} \red{betti\_numbers}(simplicial\_complex):\\
							\text{\quad\quad}{\color{gray}{\# Create a directed graph from the simplicial complex}}\\
							\text{\quad\quad}G = nx.DiGraph()\\
							\\
							\text{\quad\quad}\blue{for} simplex \blue{in} simplicial\_complex:\\
							\text{\quad\quad}\text{\quad\quad}\blue{for} i \blue{in} {\color{orange}range}({\color{orange}len}(simplex)):\\
							\text{\quad\quad}\text{\quad\quad\quad\quad}\blue{for} j \blue{in} {\color{orange}range}(i + \red{1}, {\color{orange}len}(simplex)):\\
							\text{\quad\quad}\text{\quad\quad\quad\quad\quad\quad}G.add\_edge(simplex[i], simplex[j])\\
							\\
							\text{\quad\quad}{\color{gray}\# Compute the directed graph Laplacian matrix}\\
							\text{\quad\quad}laplacian\_matrix = nx.directed\_laplacian\_matrix(G).toarray()\\
							\\
							\text{\quad\quad}{\color{gray}\# Compute the matrix rank}\\
							\text{\quad\quad}rank = scipy.linalg.matrix\_rank(laplacian\_matrix)\\
							\text{\quad\quad}{\color{gray}{\# Compute the Betti numbers}}\\
							\text{\quad\quad}num\_nodes = {\color{orange}len}(G.nodes())\\
							\text{\quad\quad}betti\_0 = num\_nodes - rank\\
							\text{\quad\quad}betti\_1 = {\color{orange}len}(simplicial\_complex)- num\_nodes + rank\\
							\\
							\text{\quad\quad}\blue{return} betti\_0, betti\_1\\
							{\color{gray}{\# Example usage}}\\
							\blue{if} \_\_name\_\_ == {\color{green}"\_\_main\_\_"}:\\
							\text{\quad\quad}{\color{gray}\# Define your simplicial complex as a list of simplices (lists of vertices)}\\
							\text{\quad\quad}simplicial\_complex = [[\red{O}, \red{1}],[\red{1},\red{2}],[\red{2},\red{ 3}],[\red{3},\red{4}],[\red{0},\red{2},\red{4}]]\\
							\text{\quad\quad}betti\_0, betti\_1 = betti\_numbers(simplicial\_complex)
							\\
							{\color{orange}print}({\color{green}"Betti 0:"},betti\_0)\\
							{\color{orange}print}({\color{green}"Betti 1:"},betti\_1)
					}}
				\end{mdframed}
				
				In this code:\\
				1. We create a directed graph from the simplicial complex, where each node represents a vertex, and directed edges connect vertices if they are part of the same simplex.\\
				2. We compute the directed graph's Laplacian matrix using \textbf{`nx.directed\_laplacian\_matrix'}.\\
				3. We calculate the rank of the Laplacian matrix to determine the number of independent cycles in the complex.\\
				4. Finally, we use the rank information to compute the Betti numbers: Betti 0 and Betti 1.\\
				You can modify the \textbf{`simplicial\_complex'} variable to represent your specific simplicial complex, and the code will compute the corresponding Betti numbers.
			\end{mdframed}
			\caption{ChatGPT provides the answer for computing Betti numbers (I)}\label{fig:direct_ask_1}
		\end{figure}

\begin{figure}[H]
			\centering
			\begin{mdframed}[backgroundcolor=lightgray, linewidth=0pt]
				\noindent Q: Given an abstract simplicial complex, can you give me a Python code for computing Betti numbers for the simplicial complex?
			\end{mdframed}
			\begin{mdframed}[backgroundcolor=orange!30, linewidth=0pt]
				\noindent \includegraphics[width=0.15in]{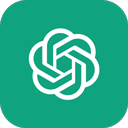} 
				Certainly! To compute the Betti numbers for an abstract simplicial complex in Python, you can use libraries such as NetworkX or Gudhi. Here, I'll show you how to compute Betti numbers using the Gudhi library, which is specifically designed for topological data analysis.\\
				First, you'll need to install Gudhi if you haven't already. You can install it using pip:				
				\begin{mdframed}[backgroundcolor=black, linewidth=0pt]
					\textcolor[rgb]{1.00,1.00,1.00}	{
						\texttt{pip install gudhi
					}}
				\end{mdframed}
				Now, let's compute the Betti numbers for your simplicial complex:
				\begin{mdframed}[backgroundcolor=black, linewidth=0pt]
					\textcolor[rgb]{1.00,1.00,1.00}	{
						\texttt{\blue{import} gudhi\\
							\gray{\# Create a simplicial complex}\\
							simplicial\_complex = gudhi.SimplexTree()\\
							\\
							\gray{\# Add simplices to the complex}\\
							\gray{\# For example, adding a triangle}\\
							simplicial\_complex.insert([\red{0}, \red{1}, \red{2}])\\
							\gray{\# Add more simplices as needed}\\
							\\
							\gray{\# Compute Betti numbers}\\
							betti\_numbers = simplicial\_complex.betti\_numbers()\\
							\\
							\gray{\# Print the Betti numbers}\\
							\orange{print}(\green{"Betti numbers:"}, betti\_numbers)
					}}
				\end{mdframed}
				Replace the \textbf{`insert'} calls with the simplices of your simplicial complex. The code above will compute the Betti numbers for the complex and print them.\\
				Here's a brief explanation of the code:\\
				1. We import the \textbf{`gudhi'} library.\\
				2. We create a \textbf{`SimplexTree'} object to represent the simplicial complex.\\
				3. We add simplices to the complex using the \textbf{`insert'} method. You should replace this with the simplices from your abstract simplicial complex.\\
				4. We compute the Betti numbers using the \textbf{`betti\_numbers'} method.\\
				5. Finally, we print the Betti numbers.\\
				This code should work for computing Betti numbers of a simplicial complex using Gudhi. Make sure to replace the simplices with your actual data.
			\end{mdframed}
			\caption{ChatGPT provides the answer for computing Betti numbers (II)}\label{fig:direct_ask_2}
		\end{figure}

		In this context, our goal is to utilize  ChatGPT to provide computer codes that take an abstract simplicial complex as input and calculates its Betti numbers as output. However, if we directly present this task to ChatGPT, it might offer inappropriate or impractical coding suggestions, or it might recommend using existing libraries. The former is not what we desire, and the latter does not align with our goal of modifying and adapting the code to suit our specific needs. As shown in Fig. \ref{fig:direct_ask_1} and Fig. \ref{fig:direct_ask_2}, for the same question, ChatGPT provides different responses along with corresponding code snippets. In its responses, ChatGPT presents a Python code, but it references external libraries such as NetworkX and Gudhi. Importantly, referencing the Gudhi library includes functions for computing simplicial complexes and their Betti numbers, which is not the approach we want. Additionally, the method suggested in the response to Fig. \ref{fig:direct_ask_1}, involving the 0-dimensional Laplacian matrix of the simplicial complex, is also incorrect as it cannot help us compute the Betti numbers of the complex.
		
		The example mentioned above highlights the necessity of utilizing our own mathematical understanding and expertise to guide and supervise ChatGPT in writing mathematical codes. This involves two main aspects. First, we need to train ChatGPT to ``understand'' the concepts and notations we provide, such as simplicial complexes, homology, and Betti numbers. Second, we need to design test problems for which we know the answer and guide ChatGPT to come up with codes to compute correct results.

		With regard to this matter, we intended to train ChatGPT to use the formula in Eq. (\ref{equ:betti}) to calculate the Betti numbers of a simplicial complex. During this process, we needed  to make sure ChatGPT understood what an abstract simplicial complex is. To achieve this, we began by asking ChatGPT for the definition of an abstract simplicial complex. We then discovered that ChatGPT can comprehend the concept of an abstract simplicial complex. However, we had noticed that the representation of simplicial complexes may vary. For example, the representation of $\Delta[1]$ could sometimes be given as $[\{0\},\{1\},\{0,1\}]$, or $[[0],[1],[0,1]]$, or $[(0),(1),(0,1)]$. And the order of simplices might not always follow lexicographical order. These variations in representation might be due to differences in computer data types or mathematical expressions.		
		To address this issue, we had decided to standardize the representation of simplicial complexes to a format similar to $[[0],[1],[0,1]]$, with simplices ordered in lexicographical order. We did not dwell on the specific computer data types and storage mechanisms.
		
		
		After confirming that ChatGPT can correctly understand abstract simplicial complexes, we decided to guide it in generating abstract simplicial complexes in a specific manner. More precisely, by providing a finite number of points in a Euclidean space and a real-value distance as input, ChatGPT needed to generate a Rips complex as output.
		This process may not succeed in a single attempt, but through iterative questioning and some debugging, ChatGPT successfully provided the desire code. The sample results provided by ChatGPT's code is shown in Fig. \ref{fig:show_1}. In this example, we take the discrete point set $\{(0, 0), (2, 1), (2, 1), (3, 3), (3 ,4)\}$ and the distance $d=3.0$ as input.
		\begin{figure}[H]
			\centering
			\begin{mdframed}[backgroundcolor=blue!30, linewidth=0pt]
				\noindent   Abstract Simplicial Complex:\\
				Vertices: [0, 1, 2, 3, 4]\\
				Simplices: [[0], [1], [2], [3], [4], [0, 1], [0, 2], [1, 2], [1, 3], [2, 3], [3, 4], [0, 1, 2], [1, 2, 3]]\\
				Highest Dimension: 2
			\end{mdframed}
			\caption{An example for generating a Rips complex}\label{fig:show_1}
		\end{figure}
		
\begin{figure}[H]
			\centering
			\begin{mdframed}[backgroundcolor=lightgray, linewidth=0pt]
				\noindent Q: Given an abstract simplicial complex [[2], [5], [4], [1], [3], [2, 4], [1, 2], [3, 4], [2, 3], [2, 3, 4]], can you show how to compute the boundary matrix by a Python code?
			\end{mdframed}
			\begin{mdframed}[backgroundcolor=orange!30, linewidth=0pt]
				\noindent \includegraphics[width=0.15in]{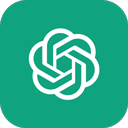} 
				To compute the boundary matrix of an abstract simplicial complex, you can use Python. The boundary matrix represents the boundary operators that map simplices to their faces. Here's a Python code to compute the boundary matrix for the given abstract simplicial complex:
				\begin{mdframed}[backgroundcolor=black, linewidth=0pt]
					\textcolor[rgb]{1.00,1.00,1.00}	{
						\texttt{\blue{import} numpy \blue{as} np\\
							\\
							\gray{\# Define the abstract simplicial complex as a list of simplices}\\
							complex = [[\red{2}], [\red{5}], [\red{4}], [\red{1}], [\red{3}], [\red{2}, \red{4}], [\red{1}, \red{2}], [\red{3}, \red{4}], [\red{2}, \red{3}], [\red{2}, \red{3}, \red{4}]]\\
							\\
							\gray{\# Extract the unique vertices from the simplicial complex}\\
							vertices = \orange{list}(\orange{set}([v \blue{for} simplex \blue{in} complex \blue{for} v \blue{in} simplex]))\\
							vertices.sort() \gray{\# Sort vertices for consistent ordering}\\
							\\
							\gray{\# Create a dictionary to map vertices to their indices}\\
							vertex\_to\_index = \{v: i \blue{for} i, v \blue{in} \orange{enumerate}(vertices)\}\\
							\\
							\gray{\# Determine the number of vertices and simplices}\\
							num\_vertices = \orange{len}(vertices)\\
							num\_simplices = \orange{len}(complex)\\
							\\
							\gray{\# Initialize the boundary matrix with zeros}\\
							boundary\_matrix = np.zeros((num\_simplices, num\_vertices), dtype=\orange{int})\\
							\\
							\gray{\# Fill in the boundary matrix}\\
							\blue{for} i, simplex \blue{in} \orange{enumerate}(complex):\\
							\text{\quad\quad}\blue{for} j, vertex \blue{in} \orange{enumerate}(simplex):\\
							\text{\quad\quad}\text{\quad\quad}\gray{\# Remove the j-th vertex from the simplex to get the (j-1)-face}\\
							\text{\quad\quad}\text{\quad\quad}face = simplex[:j] + simplex[j + \red{1}:]\\
							\\
							\text{\quad\quad}\text{\quad\quad}\gray{\# Find the index of the face in the vertex list}\\
							\text{\quad\quad}\text{\quad\quad}face\_index = vertex\_to\_index[face[\red{0}]]\\
							\\
							\text{\quad\quad}\text{\quad\quad}\gray{\# Add 1 to the (i, face\_index) entry of the boundary matrix}\\
							\text{\quad\quad}\text{\quad\quad}boundary\_matrix[i, face\_index] = \red{1} \blue{if} j \% \red{2} == \red{0} \blue{else} -\red{1}\\
							\\
							\gray{\# Print the boundary matrix}\\
							\orange{print}(\green{"Boundary Matrix:"})\\
							\orange{print}(boundary\_matrix)
					}}
				\end{mdframed}
				This code first extracts the unique vertices and assigns indices to them, initializes the boundary matrix with zeros, and then fills in the matrix according to the boundary operator rules. The result is the boundary matrix for the given abstract simplicial complex.
			\end{mdframed}
			\caption{We ask ChatGPT how to compute the boundary matrix.}\label{fig:boundary}
		\end{figure}

		In order to enable ChatGPT to generate codes for computing Betti numbers correctly, we had decided to have it first performing the computation of the boundary matrix. As before, we first confirmed whether ChatGPT understood what the boundary operator on the chain complex associated with a simplicial complex is. Fig. \ref{fig:boundary} shows how we approached ChatGPT with this task.
		
		ChatGPT provided us with various responses to the task we were given. As object-oriented individuals, we chose our preferences to proceed with the question. The result code from ChatGPT's responses is shown in Fig. \ref{fig:boundary}. We had placed and executed the code into Python, however, obtaining an error as shown in Fig. \ref{fig:error_1}.
		\begin{figure}[H]
			\centering
			\begin{mdframed}[backgroundcolor=blue!30, linewidth=0pt]
				\noindent   File ``C:\textbackslash Users\textbackslash Lenovo\.spyder-py3\textbackslash ChatGPT for topology\textbackslash boundary\_2.py'', line 35, in $<$module$>$\\
				face\_index = vertex\_to\_index[face[0]]\\
				IndexError: list index out of range
			\end{mdframed}
			\caption{Feedback information on the error.}\label{fig:error_1}
		\end{figure}
		For the debugging process, we directly input the error message to ChatGPT for feedback. After several attempts, ChatGPT offered the code without errors by adding the code in Fig. \ref{fig:debug_1} before ``\# Find the index of the face in the vertex list'' in Fig. \ref{fig:boundary}.
		\begin{figure}[H]
			\centering
			\begin{mdframed}[backgroundcolor=black, linewidth=0pt]
				\textcolor[rgb]{1.00,1.00,1.00}{ \texttt{\blue{if} \orange{len}(face) == \red{0}:\\
						\text{\hspace{1em}} \blue{continue}}}
			\end{mdframed}
			\caption{The code added after debugging.}\label{fig:debug_1}
		\end{figure}
		
		We validated the code through several test examples and found that the computed boundary matrices were incorrect. At this point, we need to review the code and understand what caused the error. Fortunately, ChatGPT provides comments for its codes, which will help us read and understand it better. Additionally, through interaction with ChatGPT, we can inquire about the parts of the code that are difficult to understand, and ChatGPT will respond accordingly. This process can not only help us quickly pinpoint the reasons for code errors but also evaluate whether ChatGPT lacks  the understanding of mathematical algorithms.
		ChatGPT plays the role of providing a general code framework and interpreting the code in the above process. In the end, we obtained the code as shown in Fig. \ref{fig:betti}, which can correctly compute Betti numbers of a given abstract simplicial complex.

		\begin{figure}[H]
			\centering
			\begin{mdframed}[backgroundcolor=black, linewidth=0pt]
				\textcolor[rgb]{1.00,1.00,1.00}	{					
					\texttt{\blue{import} numpy \blue{as} np\\
						\\
						\blue{def} \red{compute\_betti}(complex):\\
						\text{\quad\quad}num\_simplices = {\color{orange}len}(complex)\\
						\text{\quad\quad}m = {\color{orange}len}(complex[num\_simplices - \red{1}]) - \red{1}\\
						\text{\quad\quad}b = [\red{0}, \red{0}]\\
						\text{\quad\quad}k = \red{1} \\
						\\
						\text{\quad\quad}\blue{for} i \blue{in} {\color{orange}{range}}(num\_simplices):\\
						\text{\quad\quad}\text{\quad\quad}\blue{if} {\color{orange}len}(complex[i]) == k:\\
						\text{\quad\quad}\text{\quad\quad}\text{\quad\quad}b[k] += \red{1}\\
						\text{\quad\quad}\text{\quad\quad}{\color{orange}else}:\\
						\text{\quad\quad}\text{\quad\quad}\text{\quad\quad}k += \red{1}\\
						\text{\quad\quad}\text{\quad\quad}\text{\quad\quad}b.append(b[k - \red{1}] + \red{1})\\
						\\
						\text{\quad\quad}boundary\_matrix = np.zeros((num\_simplices, num\_simplices), dtype={\color{orange}int})\\
						\\
						\text{\quad\quad}\blue{for} i, simplex \blue{in} {\color{orange}enumerate}(complex):\\
						\text{\quad\quad}\text{\quad\quad}\blue{for} j, vertex \blue{in} {\color{orange}enumerate}(simplex):\\
						\text{\quad\quad}\text{\quad\quad}\text{\quad\quad}face = simplex[:j] + simplex[j + \red{1}:]\\
						\text{\quad\quad}\text{\quad\quad}\text{\quad\quad}face\_index = complex.index(face)\\
						\text{\quad\quad}\text{\quad\quad}\text{\quad\quad}boundary\_matrix[i, face\_index] = \red{1} \blue{if} j \% \red{2} == \red{0} else -\red{1}\\
						\\
						\text{\quad\quad}r = [\red{0}] * (m + \red{1})\\
						\text{\quad\quad}\blue{for} i \blue{in} {\color{orange}{range}}(m):\\
						\text{\quad\quad}\text{\quad\quad}r[i + \red{1}] = np.linalg.matrix\_rank(boundary\_matrix[b[i + \red{1}]:b[i + \red{2}], b[i]:b[i + 1]])\\
						\\
						\text{\quad\quad}betti = [\red{0}] * (m + \red{1})\\
           			    \text{\quad\quad}\blue{if} m == \red{0}:\\
                        \text{\quad\quad\quad\quad}betti[\red{0}] = b[\red{1}] - b[\red{0}]\\
                        \text{\quad\quad}\blue{elif} m == \red{1}:\\
                        \text{\quad\quad\quad\quad}betti[\red{0}] = b[\red{1}] - b[\red{0}] - r[\red{0}] - r[\red{1}]\\
                        \text{\quad\quad\quad\quad}betti[m] = b[\red{2}] - b[\red{1}] - r[\red{1}]\\
                        \text{\quad\quad}\blue{else}:\\
                        \text{\quad\quad\quad\quad} \blue{for} i \blue{in} {\color{orange}{range}}(m):\\
                        \text{\quad\quad\quad\quad\quad}\text{\quad\quad}betti[i] = b[i + \red{1}] - b[i] - r[i] - r[i + \red{1}]\\
						\\
						\text{\quad\quad}\blue{return} betti\\
						\\
						{\color{gray}\# Define the abstract simplicial complex as a list of simplices}\\
						complex = [[\red{1}], [\red{2}], [\red{3}], [\red{4}], [\red{5}], [\red{1}, \red{2}], [\red{2}, \red{3}], [\red{2}, \red{4}], [\red{3}, \red{4}], [\red{2}, \red{3}, \red{4}]]\\
						\\
						{\color{gray}{\# Calculate Betti numbers using the function}}\\
						betti\_numbers = compute\_betti(complex)\\
						\\
						{\color{gray}\# Print the Betti numbers}\\
						{\color{orange}print}({\color{green}"Betti Numbers:"})\\
						{\color{orange}print}(betti\_numbers)
				}}
			\end{mdframed}
			\caption{The final code for computing Betti numbers.}\label{fig:betti}
		\end{figure}
		
		\subsection{Computing Dirac and Laplacian matrices}
		
		In this section, we further guide ChatGPT to generate the codes for computing Dirac and Laplacian operators of a simplicial complex.
		
		Let $\mathcal{K}$ be a simplicial complex, and let $C_{\ast}(\mathcal{K})$ be chain complex of $\mathcal{K}$. One can endow $C_{\ast}(\mathcal{K})$ with an inner product structure given by
		\begin{equation*}
			\langle\sigma,\tau\rangle=\left\{
			\begin{array}{ll}
				1, & \hbox{$\sigma=\tau$;} \\
				0, & \hbox{otherwise.}
			\end{array}
			\right.
		\end{equation*}
		Then we have the adjoint operator $d^{\ast}:C_{\ast}(\mathcal{K})\to C_{\ast}(\mathcal{K})$ of $d$ with respect to the above inner product. The \emph{Dirac operator} $\delta:C_{\ast}(\mathcal{K})\to C_{\ast}(\mathcal{K})$ on $\mathcal{K}$ is defined by $\delta=d+d^{\ast}$. The \emph{Laplacian operator} $\Delta:C_{\ast}(\mathcal{K})\to C_{\ast}(\mathcal{K})$ on $\mathcal{K}$ is defined by $\Delta=\delta^{2}=d^{\ast}d+dd^{\ast}$. As we know, the eigenvalues of Laplacian operators are crucial features for capturing the shape of a data-set. The number of zero eigenvalues is equal to the Betti number, and the second smallest eigenvalue, called the Fiedler vector, corresponds to the connectivity of the simplicial complex. The smallest positive eigenvalue, denoted by $\lambda(1)$ is the spectral gap, which is closely related the Cheeger constant in differential geometry. The spectral gap $\lambda(1)$ is also one of the most commonly used features in applications.
		
		Let $B_{p}$ be the representation matrix of $d_{p}:C_{p}(\mathcal{K})\to C_{p-1}(\mathcal{K})$. More precisely, we write
        \begin{equation*}
          d_{p}\left(
                      \begin{array}{c}
                        e_{p}^{1} \\
                        \vdots \\
                        e_{p}^{k_{p}} \\
                      \end{array}
                    \right)=B_{p}\left(
                      \begin{array}{c}
                        e_{p-1}^{1} \\
                        \vdots \\
                        e_{p-1}^{k_{p-1}} \\
                      \end{array}
                    \right)
        \end{equation*}
        Here, $e_{p}^{1},\dots,e_{p}^{k_{p}}$ is the basis of $C_{p}(\mathcal{K})$. Then the representation matrix of the $p$-dimensional Laplacian operator $\Delta_{p}:C_{p}(\mathcal{K})\to C_{p}(\mathcal{K})$ can be written as		
		\begin{equation*}
        L_{p}=B_{p}B_{p}^{T}+B_{p+1}^{T}B_{p+1},\quad p\geq 1.
        \end{equation*}
        In particular, $L_{0}=B_{1}B_{1}^{T}$.
It is worth noting that if the represention matrix $B_{p}$ for $d_{p}$ is given by right multiplication, the corresponding Laplacian matrix is $L_{p}=B_{p}^{T}B_{p}+B_{p+1}B_{p+1}^{T}$.

\begin{figure}[H]
			\centering
			\begin{mdframed}[backgroundcolor=lightgray, linewidth=0pt]
				\noindent Q: The Dirac matrix is given by the sum of the boundary matrix and its transpose matrix. Can you provide the function for computing the Dirac matrix in coding?
			\end{mdframed}
			\begin{mdframed}[backgroundcolor=orange!30, linewidth=0pt]
				\noindent \includegraphics[width=0.15in]{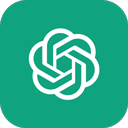} 
				Certainly! You can compute the Dirac matrix by summing the boundary matrix and its transpose. Here's a Python function to do that:
				\begin{mdframed}[backgroundcolor=black, linewidth=0pt]
					\textcolor[rgb]{1.00,1.00,1.00}	{
						\texttt{$\cdots\cdots$\\
							\blue{def} \red{compute\_dirac\_matrix}(complex):\\
							\text{\qquad}\gray{\# Compute the boundary matrix}\\
							\text{\qquad}boundary\_matrix = compute\_boundary\_matrix(complex)\\
							\text{\qquad}\gray{\# Compute the Dirac matrix by adding the boundary matrix and its transpose}\\
							\text{\qquad}dirac\_matrix = boundary\_matrix + boundary\_matrix.T\\
							\text{\qquad}\blue{return} dirac\_matrix\\
							$\cdots\cdots$
					}}
				\end{mdframed}
				This code defines a \textbf{`compute\_dirac\_matrix'} function that takes your abstract simplicial complex as input, computes the boundary matrix using the previously defined \textbf{`compute\_boundary\_matrix'} function, and then computes the Dirac matrix by adding the boundary matrix to its transpose.
			\end{mdframed}
			\caption{The added code for Dirac matrix.}\label{fig:dirac}
		\end{figure}

\begin{figure}[H]
			\centering
			\begin{mdframed}[backgroundcolor=lightgray, linewidth=0pt]
				\noindent Q: The Laplacian matrix is given by the Dirac matrix multiplies its transpose. Could you  please provide the Laplacian matrix in coding?
			\end{mdframed}
			\begin{mdframed}[backgroundcolor=orange!30, linewidth=0pt]
				\noindent \includegraphics[width=0.15in]{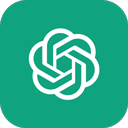} 
				To compute the Laplacian matrix, you can multiply the Dirac matrix by its transpose. Here's a Python function to do that:
				\begin{mdframed}[backgroundcolor=black, linewidth=0pt]
					\textcolor[rgb]{1.00,1.00,1.00}	{
						\texttt{$\cdots\cdots$\\
							\blue{def} \red{compute\_laplacian\_matrix}(complex):\\
							\text{\qquad}\gray{\# Compute the Dirac matrix}\\
							\text{\qquad}dirac\_matrix = compute\_dirac\_matrix(complex)\\
							\text{\qquad}\gray{\# Compute the Laplacian matrix by multiplying the Dirac matrix by its transpose}\\
							\text{\qquad}laplacian\_matrix = np.dot(dirac\_matrix, dirac\_matrix.T)\\
							\text{\qquad}\blue{return} laplacian\_matrix\\
							$\cdots\cdots$
					}}
				\end{mdframed}
				This code defines a \textbf{`compute\_laplacian\_matrix'} function that takes your abstract simplicial complex as input, computes the Dirac matrix using the previously defined \textbf{`compute\_dirac\_matrix'} function, and then computes the Laplacian matrix by multiplying the Dirac matrix by its transpose.
			\end{mdframed}
			\caption{The added code for Laplacian matrix.}\label{fig:laplacian}
		\end{figure}

		After we provided the matrix representation of Dirac and Laplacian operators to ChatGPT, it immediately generates the desire code for computing the corresponding Dirac and Laplacian matrices. The interaction with ChatGPT is shown in Fig \ref{fig:dirac} and Fig. \ref{fig:laplacian}. Although ChatGPT has been shown to perform relatively well in complex mathematical computations, its understanding of mathematical formulas is satisfactory.

		\section{ChatGPT for hypergraphs and directed graphs}\label{section:hypergraph_laplacian}
		
		Topological hypergraph and digraph (directed graph) are relatively new mathematical models for studying the complex systems in applications. From a combinatorial perspective, hypergraphs are an extension of graphs where edges can be of higher dimension. From a geometric standpoint, hypergraphs are a generalization of abstract simplicial complexes, relaxing the condition of face operator closure that is inherent in simplicial complexes. In \cite{bressan2019embedded}, the authors introduced the embedded homology hypergraphs, which are exactly a generalization of the simplicial homology of simplicial complexes. In this section, the topology invariants for hypergraphs is based on the embedded homology. The GLMY theory was introduced to study the topology of directed graphs \cite{grigor2012homologies,grigor2014homotopy,grigor2015cohomology}. Specifically, simplicial complexes are essential path complexes on digraphs, and the path homology can also be viewed as a generalization of simplicial homology. The embedded homology and GLMY theory provide us with a broader perspective and general tools to address practical problems.

		\subsection{Betti numbers for hypergraphs}\label{section:hyper_betti}
	
		Let us recall the embedded homology of hypergraphs. Let $V$ be a finite and ordered nonempty set. Let $\mathbf{P}(V)$ be the power set of $V$. A \emph{hypergraph} on $V$ is a pair $(V,\mathcal{H})$ such that $\mathcal{H}$ is a nonempty subset of $\mathbf{P}(V)$. The \emph{simplicial closure} is defined by $\Delta \mathcal{H}=\{\sigma|\sigma\subseteq \tau\in \mathcal{H}\}$, which can be viewed as the minimal simplicial complex contain $\mathcal{H}$.
		
		Let $\mathbb{K}$ be the ground field. Let $C_{\ast}(\Delta\mathcal{H};\mathbb{K})$ be the chain complex of the simplicial complex $\Delta \mathcal{H}$. Let $D_{\ast}(\mathcal{H};\mathbb{K})$ be the graded subspace of $C_{\ast}(\Delta\mathcal{H};\mathbb{K})$ generated by the hyperedges in $\mathcal{H}$. Then one has the \emph{infimum chain complex}
		\begin{equation*}
			\Inf_{\ast}(\mathcal{H};\mathbb{K})=\sum\{C'_{\ast}\subseteq D_{\ast}(\mathcal{H};\mathbb{K})|C'_{\ast}\text{ is a subchain complex of }C_{\ast}(\Delta\mathcal{H};\mathbb{K})\}
		\end{equation*}
		and the \emph{supremum chain complex}
		\begin{equation*}
			\Sup_{\ast}(\mathcal{H};\mathbb{K})=\bigcap\{C'_{\ast}\supseteq D_{\ast}(\mathcal{H};\mathbb{K})|C'_{\ast}\text{ is a subchain complex of }C_{\ast}(\Delta\mathcal{H};\mathbb{K})\}
		\end{equation*}
		of $\mathcal{H}$. It is proved that the homology of the infimum chain complex is isomorphic to the supremum chain complex of $\mathcal{H}$. The \emph{embedded homology} is defined to be
		\begin{equation*}
			H_{p}(\mathcal{H};\mathbb{K}):=H_{p}(\Inf_{\ast}(\mathcal{H};\mathbb{K})),\quad p\geq 0.
		\end{equation*}

\begin{figure}[H]
			\centering
			\begin{mdframed}[backgroundcolor=lightgray, linewidth=0pt]
				\noindent Q: Find the indices of each element in a subset within another ordered set's subset and provide Python code for it.
			\end{mdframed}
			\begin{mdframed}[backgroundcolor=orange!30, linewidth=0pt]
				\noindent \includegraphics[width=0.15in]{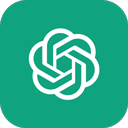} 
				To find the indices of each element in a subset within another ordered set's subset in Python, you can use list comprehensions and the \textbf{`enumerate()'} function. Here's a Python code example to achieve this:
				\begin{mdframed}[backgroundcolor=black, linewidth=0pt]
					\textcolor[rgb]{1.00,1.00,1.00}	{\texttt{\blue{def} \red{find\_indices\_of\_subset}(main\_set, subset):\\
							\text{\quad\quad}indices = [i \blue{for} i, x in \orange{enumerate}(main\_set) \blue{if} x \blue{in} subset]\\
							\text{\quad\quad}\blue{return} indices\\
							\\
							\gray{\# Example usage:}\\
							main\_set = [\red{1}, \red{2}, \red{3}, \red{4}, \red{5}, \red{6}, \red{7}, \red{8}, \red{9}]\\
							subset = [\red{3}, \red{5}, \red{7}]\\
							\\
							indices = find\_indices\_of\_subset(main\_set, subset)\\
							\orange{print}(\green{"Indices of subset elements in the main set:"}, indices)}}
				\end{mdframed}
				In this code:\\
				1. We define a function \textbf{`find\_indices\_of\_subset'} that takes two lists as input, \textbf{`main\_set'} and \textbf{`subset'}.\\
				2. We use a list comprehension to iterate over the elements in the \textbf{`main\_set'} along with their indices using the \textbf{`enumerate()'} function.\\
				3. For each element \textbf{`x'} in \textbf{`main\_set'}, we check if it exists in the \textbf{`subset'}. If it does, we append its index \textbf{`i'} to the \textbf{`indices'} list.\\
				4. Finally, we return the \textbf{`indices'} list, which contains the indices of elements in the \textbf{`subset'} that appear in the \textbf{`main\_set'}.
			\end{mdframed}
			\caption{ChatGPT helps us find the index of a subset within an ordered set.}\label{fig:find_index}
		\end{figure}

		Now, we will sketch the algorithm for computing the Betti numbers for hypergraphs. By \cite[Proposition 2.1]{bressan2019embedded}, one has that $\Inf_{\ast}(\mathcal{H};\mathbb{K})=D_{\ast}(\mathcal{H};\mathbb{K})\cap d^{-1}D_{\ast}(\mathcal{H};\mathbb{K})$. Here, $d$ is the differential on $C_{\ast}(\Delta\mathcal{H};\mathbb{K})$. In this work, the calculation of embedded homology is based on the chain complex $D_{\ast}(\mathcal{H};\mathbb{K})\cap d^{-1}D_{\ast}(\mathcal{H};\mathbb{K})$. By \cite[Theorem 2.3]{liu2022neighborhood}, we have
		\begin{equation*}
			\beta_{p}(\mathcal{H})=\dim D_{p}(\mathcal{H};\mathbb{K})\cup dD_{p+1}(\mathcal{H};\mathbb{K})-\dim dD_{p+1}(\mathcal{H};\mathbb{K}) -\dim dD_{p}(\mathcal{H};\mathbb{K}).
		\end{equation*}
		Let $B_{p}$ be the representation matrix of the differential $d_{p}:D_{p}(\mathcal{H};\mathbb{K})\to C_{p-1}(\Delta\mathcal{H};\mathbb{K})$. Then the Betti numbers can be computed by
		\begin{equation}\label{equ:hypergraph}
			\beta_{p}(\mathcal{H})=\rank\left(
			\begin{array}{c}
				E_{p} \\
				B_{p} \\
			\end{array}
			\right)
			-\rank B_{p+1}-\rank B_{p}.
		\end{equation}
		Here, $E_{p}$ is the matrix of $p$-hyperedges. We need to compute the boundary matrix $B_{p}$ and the hyperedge matrix $E_{p}$. The algorithm for computing the boundary matrix of a hypergraph is almost identical to the algorithm for computing the boundary matrix of a simplicial complex. The key difference lies in the induces of hyperedges within the simplicial closure of the hypergraph. More precisely, we should find the index of each hyperedge in the list of the simplices of the simplicial complex $\Delta \mathcal{H}$. Besides, we also need the index of the boundary of each hyperedge. Based on the reasons mentioned above, we posed the question to ChatGPT in Fig. \ref{fig:find_index}.
		
		As shown Fig. \ref{fig:dirac} and Fig. \ref{fig:laplacian}, we can obtain the code for computing the Betti numbers for hypergraphs based on Eq. (\ref{equ:hypergraph}). In the progress, we need to   validate whether the examples computed by the algorithm are correct. We examined the computed results at various locations in the code and determined which step of the calculation went wrong. We also tried to understand what caused these errors. While ChatGPT plays a significant role in writing and explaining code, our primary expertise lies in comprehending mathematical principles, and our responsibility is to oversee ChatGPT rather than engage in coding.

		\subsection{Betti numbers for directed graphs}\label{section:Betti_digraph}
	
		Directed graphs (digraphs) naturally possess a 1-complex structure but lack higher-dimensional structure. However, the GLMY theory can detect the higher-dimensional structure of digraphs and compute their higher-order homology, referred to as path homology. Persistent path homology has been applied in various field \cite{chen2023path,chowdhury2022path,chowdhury2018persistent}.

From now on, the Betti numbers for digraphs always refer to the Betti numbers based on the path homology.
		The algorithm for computing the Betti numbers for digraphs is quite similar. Recall that a digraph $G=(V,E)$ consists of a finite nonempty set $V$, called the vertex set, and the edge set $E\subseteq V\times V$. From now on, all the digraphs considered are assumed to be simple, that is, there is no loop or multi-edge in the same direction.
		An \emph{elementary $p$-path} on $V$ is a sequence $i_{0}i_{1}\cdots i_{p}$ of vertex in $V$. Let $\Lambda_{p}(V;\mathbb{K})$ be the $K$-linear space generated by all the elementary $p$-paths on $V$. With special attention, the basis for $\Lambda(V;\mathbb{K})$ is denoted by $\{e_{i_{0}i_{1}\cdots i_{p}}\}$, where $i_{0}i_{1}\cdots i_{p}$ runs across all the elementary $p$-paths on $V$.
		It follows that $\Lambda_{\ast}(V;\mathbb{K})$ is a chain complex with the differential $d_{p}:\Lambda_{p}(V;\mathbb{K})\to \Lambda_{p-1}(V;\mathbb{K})$ given by
		\begin{equation*}
			d_{p}e_{i_{0}i_{1}\cdots i_{p}}=\sum\limits_{t=0}^{p}(-1)^{t}e_{i_{0}\cdots \hat{i_{t}}\cdots i_{p}},\quad p\geq 1
		\end{equation*}
		and $d_{0}e_{i_{0}}=0$. Here, $\hat{i_{t}}$ means omission of the index $i_{t}$. An \emph{allowed $p$-path} on a digraph $G$ is an elementary path $i_{0}i_{1}\cdots i_{p}$ such that $(i_{t-1},i_{t})\in E$ for $t=1,\dots,p$. Let $\mathcal{A}_{p}(G;\mathbb{K})$ be the $\mathbb{K}$-linear space generated by all the allowed $p$-paths on $G$. Then $\mathcal{A}_{p}(G;\mathbb{K})$ is a subspace of $\Lambda(V;\mathbb{K})$. Let
		\begin{equation*}
			\Omega_{p}(G;\mathbb{K})=\{x\in \mathcal{A}_{p}(G;\mathbb{K})|d_{p}x\in \mathcal{A}_{p-1}(G;\mathbb{K})\}.
		\end{equation*}
		Then $\Omega_{\ast}(G;\mathbb{K})$ is a chain complex. The \emph{path homology} of $G$ is defined by
		\begin{equation*}
			H_{p}(G;\mathbb{K}):=H_{p}(\Omega_{\ast}(G;\mathbb{K})),\quad p\geq 0.
		\end{equation*}
		
		It is worth noting that the chain complex $\Omega_{\ast}(G;\mathbb{K})$ is exactly the infimum chain complex of $\mathcal{A}_{\ast}(G;\mathbb{K})$ in $\Lambda_{\ast}(V;\mathbb{K})$. The construction of $\Omega_{\ast}(G;\mathbb{K})$ for digraphs is quite similar as the construction of $\Inf_{\ast}(\mathcal{H};\mathbb{K})$ for hypergraphs. This similarity determines that we can use algorithms for calculating Betti numbers of hypergraphs to compute the Betti numbers of directed graphs in GLMY theory. The main difference in computing the Betti numbers of a directed graph is that we need to find all paths and represent these paths as lists, which is consistent with the representation of simplices in simplicial complexes or hyperedges in hypergraphs. It is worth noting that there may be simplified algorithms for applying the boundary operator to paths in view of the specific characteristics of paths.
		
		Here, let us briefly outline the logic of how we will approach by asking ChatGPT to demonstrate how it can gradually help us write the code we need: Given a directed graph, provide a list of paths.
		\begin{enumerate}
			\item[(1)] Question: Can you give me a Python code for finding all the directed path on a digraph?\\
			ChatGPT gives us the all the paths for a given starting node and ending node.
			\item[(2)] Question: I need to obtain all paths, regardless of the starting and ending nodes.\\
			When we run the code, it does not show the 0-dimensional paths.
			\item[(3)] Question: Please add vertices as path of length 0 into the 'paths'.\\
			Then ChatGPT provides the satisfactory code.
		\end{enumerate}
		Besides, we also  request a code for finding the path with the  length not larger than 2 or an integer number. Indeed, a finite digraph, for example, $V=\{0,1\},E=\{(0,1),(1,0)\}$, there are infinite many paths which behave as the binary numbers.
	
		\subsection{Computing hypergraph Laplacian}\label{section:hyper_Laplacian}
	
		In this section, we will study the algorithms for computing the Dirac and Laplacian matrices of on hypergraphs and digraphs. We will present the mathematical progress for hypergraphs first. The algorithm for digraphs is much similar.
		
		Let $\mathcal{H}$ be a hypergraph. Recall that we have the chain complex $C_{\ast}(\Delta\mathcal{H};\mathbb{K})$ and its subspace $D_{\ast}(\mathcal{H};\mathbb{K})$. Let $e_{p}^{1},\dots,e_{p}^{l_{p}}$ and $\varepsilon_{p}^{1},\dots,\varepsilon_{p}^{k_{p}}$ be the basis of $C_{\ast}(\Delta\mathcal{H};\mathbb{K})$ and $D_{\ast}(\mathcal{H};\mathbb{K})$, respectively. We write the differential $d_{p}:D_{p}(\mathcal{H};\mathbb{K})\to C_{p-1}(\Delta\mathcal{H};\mathbb{K})$ by
		\begin{equation*}
			d_{p}\left(
			\begin{array}{c}
				\varepsilon_{p}^{1} \\
				\vdots\\
				\varepsilon_{p}^{k_{p}} \\
			\end{array}
			\right)=B_{p}\left(
			\begin{array}{c}
				e_{p}^{1} \\
				\vdots \\
				e_{p}^{l_{p}} \\
			\end{array}
			\right).
		\end{equation*}
		Let $\pi_{p}:C_{p}(\Delta\mathcal{H};\mathbb{K})\to C_{p}(\Delta\mathcal{H};\mathbb{K})/D_{p}(\mathcal{H};\mathbb{K})$ and $\rho_{p}:C_{p}(\Delta\mathcal{H};\mathbb{K})\to D_{p}(\mathcal{H};\mathbb{K})$ be projections of $\mathbb{K}$-linear spaces. Then we have $\Inf_{p}(\mathcal{H};\mathbb{K})=\ker (\pi\circ d_{p})$. Let $\bar{B}_{p}$ be the representation matrix of $\pi\circ d_{p}$. Let $A_{p}\left(
		\begin{array}{c}
			\varepsilon_{p}^{1} \\
			\vdots\\
			\varepsilon_{p}^{k_{p}} \\
		\end{array}
		\right)$ be an orthogonal basis of $\Inf_{p}(\mathcal{H};\mathbb{K})$. The we have
		\begin{equation*}
			\pi \circ d_{p}\left(A_{p}\left(
			\begin{array}{c}
				\varepsilon_{p}^{1} \\
				\vdots\\
				\varepsilon_{p}^{k_{p}} \\
			\end{array}
			\right)\right)=A_{p}\left(\pi \circ d_{p}\left(
			\begin{array}{c}
				\varepsilon_{p}^{1} \\
				\vdots\\
				\varepsilon_{p}^{k_{p}} \\
			\end{array}
			\right)\right)=A_{p}\bar{B}_{p}\left(
			\begin{array}{c}
				e_{p}^{1} \\
				\vdots \\
				e_{p}^{l_{p}} \\
			\end{array}
			\right).
		\end{equation*}
		Note that $A_{p}$ is an $r_{p}\times k_{p}$ matrix. Here, $k_{p}$ is the dimension of $\Inf_{p}(\mathcal{H};\mathbb{K})$. The algorithm to obtain the matrix $A_{p}$ is to calculate the left zero matrix of $\bar{B}_{p}$. Let $\tilde{B}_{p}$ be the representation matrix of $\rho\circ d_{p}$. Based on the discussion above, we have
		\begin{equation*}
			d_{p}\left(A_{p}\left(
			\begin{array}{c}
				\varepsilon_{p}^{1} \\
				\vdots\\
				\varepsilon_{p}^{k_{p}} \\
			\end{array}
			\right)\right)=\rho\circ d_{p}\left(A_{p}\left(
			\begin{array}{c}
				\varepsilon_{p}^{1} \\
				\vdots\\
				\varepsilon_{p}^{k_{p}} \\
			\end{array}
			\right)\right)=A_{p}\left(\tilde{B}_{p}\left(
			\begin{array}{c}
				\varepsilon_{p}^{1} \\
				\vdots\\
				\varepsilon_{p}^{k_{p}} \\
			\end{array}
			\right)\right)=A_{p}\tilde{B}_{p}A_{p-1}^{-1}\left(A_{p-1}\left(
			\begin{array}{c}
				\varepsilon_{p}^{1} \\
				\vdots\\
				\varepsilon_{p}^{k_{p}} \\
			\end{array}
			\right)\right).
		\end{equation*}
		Here, $A_{p}^{-1}$ is the generalized inverse or Moore-Penrose matrix of $A_{p}$, i.e., $A_{p}^{-1}A_{p}=I_{r_{p}}$. Then $A_{p}\tilde{B}_{p}A_{p-1}^{-1}$ is the desired boundary matrix of the infimum complex. Furthermore, the Laplacian matrix can be computed by
		\begin{equation*}
			L_{p}=A_{p}\tilde{B}_{p}A_{p-1}^{-1}(A_{p}\tilde{B}_{p}A_{p-1}^{-1})^{T}+(A_{p+1}\tilde{B}_{p+1}A_{p}^{-1})^{T}A_{p+1}\tilde{B}_{p+1}A_{p}^{-1},\quad p\geq 1.
		\end{equation*}
		In particular, $L_{0}=(A_{1}\tilde{B}_{1}A_{0}^{-1})^{T}A_{1}\tilde{B}_{1}A_{0}^{-1}$. The eigenvalues of Laplacian matrices are independent of the choices of the orthogonal for the infimum chain complex.
		
		Next, we   guide ChatGPT to  generate the code for computing the boundary matrix for the infimum chain complex $\Inf_{p}(\mathcal{H};\mathbb{K})$. The algorithm is based on the code for computing Betti numbers. Recall that we have obtained the code for computing the boundary matrix $B_{p}$ in Section \ref{section:hyper_betti}.

		Let $I$ be the list of the index of the hyperedges of $\mathcal{H}$ in the set of simplicies of $\Delta \mathcal{H}$.
		The matrix $\bar{B}$ is obtained by removing the columns indexed by $I$, while the matrix $\tilde{B}$ consists of the columns indexed by $I$.

		\begin{figure}[H]
			\centering
			\begin{mdframed}[backgroundcolor=lightgray, linewidth=0pt]
				\noindent Q: With a given a matrix and index list, please give me a Python code to obtain   two new matrices by removing the columns indexed by the list and not indexed by the list.
			\end{mdframed}
			\begin{mdframed}[backgroundcolor=orange!30, linewidth=0pt]
				\noindent \includegraphics[width=0.15in]{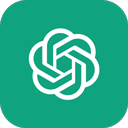} 
				To obtain two new matrices from an input matrix by removing columns based on an index list, you can use Python and NumPy. Here's a Python code snippet to achieve this:
				\begin{mdframed}[backgroundcolor=black, linewidth=0pt]
					\textcolor[rgb]{1.00,1.00,1.00}	{
						\texttt{\blue{import} numpy \blue{as} np\\
							\\
							\blue{def} \red{split\_matrix}(matrix, index\_list):\\
							\text{\qquad}\gray{\# Create an array of boolean values indicating which columns to keep}\\
							\text{\qquad}keep\_columns = np.zeros(matrix.shape[\red{1}], dtype=\orange{bool})\\
							\text{\qquad}keep\_columns[index\_list] = \blue{True}\\
							\\
							\text{\qquad}\gray{\# Use boolean indexing to split the matrix}\\
							\text{\qquad}matrix\_with\_indexed\_columns = matrix[:, keep\_columns]\\
							\text{\qquad}matrix\_without\_indexed\_columns = matrix[:, $\sim$keep\_columns]\\
							\\
							\text{\qquad}\blue{return} matrix\_with\_indexed\_columns, matrix\_without\_indexed\_columns\\
							\\
							\gray{\# Example usage:}\\
							\gray{\# Create a sample matrix and index list for testing}\\
							sample\_matrix = np.array([[\red{1}, \red{2}, \red{3}, \red{4}],[\red{5}, \red{6}, \red{7}, \red{8}],[\red{9}, \red{10}, \red{11}, \red{12}]])\\
							\\
							index\_list = [\red{0}, \red{2}]  \gray{\# Specify the indices of columns to keep}\\
							\\
							\gray{\# Call the function to split the matrix}\\
							result\_with\_indexed, result\_without\_indexed = split\_matrix(sample\_matrix, index\_list)\\
							\\
							\gray{\# Print the results}\\
							\orange{print}(\green{"Matrix with indexed columns:")}\\
							\orange{print}(result\_with\_indexed)\\
							\\
							\orange{print}(\green{"$\backslash$nMatrix without indexed columns:"})\\
							\orange{print}(result\_without\_indexed)
					}}
				\end{mdframed}
				In this code, \textbf{`split\_matrix'} takes an input matrix and an index list as arguments and returns two new matrices: one with columns indexed by the list and another with columns not indexed by the list. It uses NumPy to perform the column selection efficiently.
			\end{mdframed}
			\caption{ChatGPT helps us split a matrix into two matrices.}
		\end{figure}
		
\begin{figure}[H]
			\centering
			\begin{mdframed}[backgroundcolor=lightgray, linewidth=0pt]
				\noindent Q: Given a matrix $B$, please find a matrix $A$ such that $AB=0$ and the rows of $A$ are orthogonal to each other. Please give me a Python code.
			\end{mdframed}
			\begin{mdframed}[backgroundcolor=orange!30, linewidth=0pt]
				\noindent \includegraphics[width=0.15in]{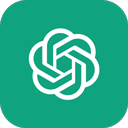} 
				$\cdots\cdots$ (Not desired content is omitted here.)
			\end{mdframed}
			\begin{mdframed}[backgroundcolor=lightgray, linewidth=0pt]
				\noindent Q: Suppose B is an $m\times n$ matrix, where $m\leq n$. Please find a $r \times m$ matrix A as above. Here, $r=m-\rank(B)$.
			\end{mdframed}
			\begin{mdframed}[backgroundcolor=orange!30, linewidth=0pt]
				\noindent \includegraphics[width=0.15in]{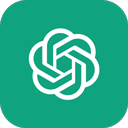} 
				If you want to find a matrix $A$ with $r$ rows (where $r = m - \rank(B)$) such that $AB = 0$ and the rows of $A$ are orthogonal to each other, you can use the concept of the left null space of matrix $B$. The left null space contains vectors that are orthogonal to the rows of $B$. Here's a Python code to do that:
				\begin{mdframed}[backgroundcolor=black, linewidth=0pt]
					\textcolor[rgb]{1.00,1.00,1.00}	{
						\texttt{\blue{import} numpy \blue{as} np\\
							\\
							\gray{\# Define matrix B}\\
							B = np.array([[\red{1}, \red{2}, \red{3}],[\red{4}, \red{5}, \red{6}],[\red{7}, \red{8}, \red{9}]])\\
							\\
							\gray{\# Calculate the rank of matrix B}\\
							rank\_B = np.linalg.matrix\_rank(B)\\
							\\
							\gray{\# Calculate r (number of rows in A)}\\
							r = B.shape[\red{0}] - rank\_B\\
							\\
							\gray{\# Perform SVD on B to get the left singular vectors (V matrix in SVD)}\\
							U, S, VT = np.linalg.svd(B)\\
							\\
							\gray{\# Extract the left null space vectors from VT}\\
							left\_null\_space\_vectors = VT[rank\_B:].T\\
							\\
							\gray{\# Create matrix A with orthogonal rows from the left null space vectors}\\
							A = left\_null\_space\_vectors[:r]\\
							\\
							\gray{\# Verify that A is orthogonal and AB = 0}\\
							result = np.dot(A, B)\\
							\orange{print}(\green{"Matrix A:"})\\
							\orange{print}(A)\\
							\orange{print}(\green{"$\backslash$nMatrix AB:"})\\
							\orange{print}(result)}}
				\end{mdframed}
				In this code:\\
				1. We define matrix $B$ as given.\\
				\qquad 2. We calculate the rank of matrix $B$ using \textbf{`numpy.linalg.matrix\_rank'}.\\
				\qquad 3. We calculate $r$, which is $m - \rank(B)$.\\
				\qquad 4. We perform Singular Value Decomposition (SVD) on $B$ to get the left singular vectors (V matrix in SVD).\\
				\qquad 5. We extract the left null space vectors from VT (transpose of V).\\
				\qquad 6. We create matrix $A$ with $r$ rows from the left null space vectors.\\
				\qquad 7. We verify that $A$ is orthogonal by calculating \textbf{`A * A.T'}, which should be the identity matrix if $A$ is orthogonal.\\
				\qquad 8. We also calculate $AB$ to ensure that $AB = 0$.\\
				Make sure you have NumPy installed (\textbf{`pip install numpy'}) to run this code.
			\end{mdframed}
			\caption{ ChatGPT helps us find the orthogonal left zero matrix of a given matrix. Note that the code cannot run smoothly before debugging.}\label{fig:left_zero}
		\end{figure}

\begin{figure}[H]
			\centering
			\begin{mdframed}[backgroundcolor=lightgray, linewidth=0pt]
				\noindent Q: Please give me a Python code for computing the right generalized inverse matrix of a given matrix.
			\end{mdframed}
			\begin{mdframed}[backgroundcolor=orange!30, linewidth=0pt]
				\noindent \includegraphics[width=0.15in]{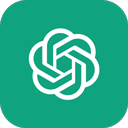} 
				To compute the right generalized inverse (also known as the Moore-Penrose inverse) of a given matrix in Python, you can use the numpy library, which provides a convenient function for this purpose. Here's how you can do it:
				\begin{mdframed}[backgroundcolor=black, linewidth=0pt]
					\textcolor[rgb]{1.00,1.00,1.00}	{
						\texttt{\blue{import} numpy \blue{as} np\\
							\\
							\gray{\# Define your matrix}\\
							A = np.array([[\red{1}, \red{2}, \red{3}],[\red{4}, \red{5}, \red{6}]])\\
							\\
							\gray{\# Compute the right generalized inverse (Moore-Penrose inverse)}\\
							A\_inv = np.linalg.pinv(A)\\
							\\
							\gray{\# Print the result}\\
							\orange{print}(\green{"Right Generalized Inverse:"})\\
							\orange{print}(A\_inv)
					}}
				\end{mdframed}
				In this code:\\
				1. We import the numpy library as \textbf{`np'}.\\
				2. We define your matrix \textbf{`A'}. You can replace this matrix with the one you want to compute the right generalized inverse for.\\
				3. We use the \textbf{`np.linalg.pinv()'} function to compute the right generalized inverse of \textbf{`A'}. This function automatically calculates the Moore-Penrose inverse.\\
				4. We print the result, which is the right generalized inverse of the input matrix \textbf{`A'}.\\
				Make sure to have numpy installed in your Python environment. You can install it using pip if it's not already installed:
				\begin{mdframed}[backgroundcolor=black, linewidth=0pt]
					\textcolor[rgb]{1.00,1.00,1.00}{ \texttt{pip install numpy}}
				\end{mdframed}
				Replace the \textbf{`A'} matrix with your desired matrix to compute its right generalized inverse.
			\end{mdframed}
			\caption{ChatGPT can provide codes for computing the generalized inverse matrix.}\label{fig:generalized_inverse}
		\end{figure}

		Then there are three tasks left for us:
		\begin{enumerate}
			\item[(1)] Compute the orthogonal left zero matrix $A_{p}$ of $\bar{B}_{p}$;
			\item[(2)] Compute the generalized inverse matrix $A_{p}^{-1}$ of $A_{p}$;
			\item[(3)] Compute the Laplacian matrices for hypergraphs.
		\end{enumerate}
		Fig. \ref{fig:left_zero} gives the code for computing orthogonal left zero matrix, while the generalized inverse matrix is provided in Fig. \ref{fig:generalized_inverse}. We combine all the functions together and obtain the final code for computing the Laplacians for hypergraphs.

		\subsection{Computing path Laplacian}

        Path Laplacian and persistent path Laplacian have also been introduced to reveal crucial structured characteristic that cannot be detected by (persistent) path homology \cite{wang2023persistent}.   Like in Section \ref{section:Betti_digraph}, the code for computing the hypergraph Laplacian is equally applicable for calculating the path Laplacian of a digraph.

        Let $G=(V,E)$ be a digraph. We have a chain complex $\Omega_{\ast}(G;\mathbb{K})$ with the differential $d_{p}:\Omega_{p}(G;\mathbb{K})\to \Omega_{p-1}(G;\mathbb{K})$. The \emph{path Laplacian} $\Delta_{p}:\Omega_{p}(G;\mathbb{K})\to \Omega_{p}(G;\mathbb{K})$ is defined by
        \begin{equation*}
          \Delta_{p}=d_{p}^{\ast}\circ d_{p}+d_{p+1}\circ d_{p+1}^{\ast},\quad p\geq 1.
        \end{equation*}
        In particular, $\Delta_{0}=d_{1}\circ d_{1}^{\ast}$. The number of zero eigenvalues of $\Delta_{p}$, known as the harmonic information, is equal to the Betti number $\beta_{p}$ from the path homology. The non-harmonic information is expressed as the positive eigenvalues of $\Delta_{p}$. Note that the smallest positive eigenvalue of $\Delta_{p}$ is one of the most used feature in illustrating the difference between path homology and path Laplacian.

        To compute the representation matrix of Laplacian operator $\Delta_{p}$, it suffices to find the representation matrix of the boundary $d_{p}$. This can be done similar as in Section \ref{section:hyper_Laplacian}. We have also obtained a functional code for computing path Laplacians.

        \subsection{Computing hyperdigraph Laplacian}

       Topological hypergraphs are highly versatile models capable of characterizing high-dimensional data structures. Topological hyperdigraphs can be viewed as assigning additional directions to each hyperedge, thereby encompassing more information and considered a more generalized mathematical model. Recently, the persistent hyperdigraph Laplacians have been introduced to characterize the topological information of complex systems \cite{chen2023persistent}.

        Let $V$ be a finite and ordered nonempty set. Let $\mathbf{S}(V)$ be the set of all sequence of distinct elements in $V$. A \emph{hyperdigraph} on $V$ is a pair $(V,\vec{\mathcal{H}})$ such that $\vec{\mathcal{H}}$ is a nonempty subset of $\mathbf{S}(V)$. A sequence in $\vec{\mathcal{H}}$ is called \emph{directed hyperedge}.
        Hyperdigraph can be viewed as a kind of generalization of hypergraphs or digraphs.

        Let $S_{p}(V;\mathbb{K})$ be the $\mathbb{K}$-linear space generated by the $(p+1)$-element sets in $\mathbf{S}(V)$. Then $S_{\ast}(V;\mathbb{K})$ is a chain complex with the boundary operator $d_{p}:S_{p}(V;\mathbb{K})\to S_{p-1}(V;\mathbb{K})$ given by
        \begin{equation*}
          d_{p}(v_{0},\dots,v_{p})=\sum\limits_{i=1}^{p}(-1)^{i}(v_{0},\dots,\hat{v_{i}},\dots,v_{p}),\quad p\geq 0.
        \end{equation*}
        Here, $\hat{v_{i}}$ means the omission of the term $v_{i}$. In particular, $d_{0}=0$ on $S_{0}(V;\mathbb{K})$. Let $F_{p}(\vec{\mathcal{H}};\mathbb{K})$ be the $\mathbb{K}$-linear space generated all the directed hyperedges in $\vec{\mathcal{H}}$.
        Let
		\begin{equation*}
			\Omega_{p}(\vec{\mathcal{H}};\mathbb{K})=\{x\in F_{p}(\vec{\mathcal{H}};\mathbb{K})|d_{p}x\in F_{p-1}(\vec{\mathcal{H}};\mathbb{K})\}.
		\end{equation*}
		Then $\Omega_{\ast}(\vec{\mathcal{H}};\mathbb{K})$ is a chain complex. The homology of $\vec{\mathcal{H}}$ is defined by
        \begin{equation*}
          H_{p}(\vec{\mathcal{H}};\mathbb{K}):=H_{p}(\Omega_{p}(\vec{\mathcal{H}};\mathbb{K})),\quad p\geq 0.
        \end{equation*}
        We endow the directed hyperedges with the standard inner product structure. The Laplacian $\Delta_{p}:\Omega_{p}(\vec{\mathcal{H}};\mathbb{K})\to \Omega_{p}(\vec{\mathcal{H}};\mathbb{K})$ is given by
        \begin{equation*}
          \Delta_{p}=d_{p}^{\ast}\circ d_{p}+d_{p+1}\circ d_{p+1}^{\ast}.
        \end{equation*}
        The representation of a hyperdigraph in the code is provided as a collection of lists, with each list representing a directed hyperedge. The Betti numbers and Laplacian matrices for hyperdigraphs can be computed in a similar manner as for hypergraphs. The ChatGPT code obtained with our supervision can also be used to compute the hyperdigraph Betti numbers and hyperdigraph Laplacians.

 \section{ChatGPT for topological persistence}\label{section:persistence}

       The development of applied topology is profoundly benefited by topological persistence. Topological persistence can characterize multi-scale topological information in data sets. Without persistence, homology theory by itself would have little application value.  Nowadays, persistent homology and persistent Laplacian are some of   the most popular TDA tools. The data set considered is always a family of discrete points in a metric space. Given a data set, one can compute the distance between any two of the points with a specific mertic, and the persistence is based on the distance. The Vietoris-Rips complex or Rips complex is the standard construction of the persistence from a data set. There have been numerous studies on persistent homology \cite{carlsson2004persistence,edelsbrunner2002topological,edelsbrunner2008persistent,zomorodian2004computing}, persistent Laplacians \cite{chen2019evolutionary,liu2023algebraic,memoli2022persistent,wang2020persistent}, and persistent Dirac operators \cite{ameneyro2022quantum,wee2023persistent}.

 \subsection{ChatGPT for persistent homology and persistent Laplacians}
        Consider the simplicial complex $\mathcal{K}_{\varepsilon}$ parametrized by a real number $\varepsilon$ such that $\mathcal{K}_{\varepsilon}\hookrightarrow \mathcal{K}_{\varepsilon'}$ for any $\varepsilon\leq \epsilon'$. For $a\leq b$, we have a morphism of chain complexes $j^{a,b}_{\ast}:C_{\ast}(\mathcal{K}_{a})\hookrightarrow C_{\ast}(\mathcal{K}_{b})$, which induces a morphism of homology groups $H_{\ast}(\mathcal{K}_{a})\to H_{\ast}(\mathcal{K}_{b})$. The \emph{$(a,b)$-persistent homology} is defined by
        \begin{equation*}
          H^{a,b}_{p}:=\im (H_{p}(\mathcal{K}_{a})\to H_{p}(\mathcal{K}_{b})),\quad p\geq 0.
        \end{equation*}
        The rank of the persistent homology $\beta_{p}^{a,b}=\rank H^{a,b}_{p}$ is the $(a,b)$-persistent Betti number, which is always represented by barcode or persistence diagram in application. In particular, we have $\beta_{p}^{a,a}=\beta_{p}^{a}$.

        The \emph{$(a,b)$-persistent Laplacian} $\Delta^{a,b}_{p}:C_{\ast}(\mathcal{K}_{a})\to C_{\ast}(\mathcal{K}_{a})$ is defined by
        \begin{equation*}
          \Delta^{a,b}_{p}=(d_{p}^{a})^{\ast}\circ d_{p}^{a}+(d_{p+1}^{a,b})^{\ast}\circ d_{p+1}^{a,b},\quad p\geq 0.
        \end{equation*}
        Here, $d_{p}^{a}$ is the boundary operator of $C_{\ast}(\mathcal{K}_{a})$, and $d_{p+1}^{a,b}$ is the composition of
        \begin{equation*}
          \{x\in C_{p+1}(\mathcal{K}_{b})|d_{p+1}^{b}x\in C_{p+1}(\mathcal{K}_{a})\}\hookrightarrow C_{p+1}(\mathcal{K}_{b})\stackrel{d_{p+1}^{b}}{\to} C_{p}(\mathcal{K}_{a}).
        \end{equation*}
        In particular, when $a=b$, the persistent Laplacian $\Delta^{a,a}_{p}$ is exactly the Laplacian $\Delta^{a}_{p}$ on chain complex $C_{\ast}(\mathcal{K}_{a})$.

        Persistent homology and persistent Laplacians can also be defined for hypergraphs, digraphs, and hyperdigraphs \cite{chen2023persistent}. This allows for the characterization of multiscale topological information in these structures.

        In applications, one often focuses on the 0th and 1st Betti numbers because they are more computationally tractable. Sometimes, researchers only compute the Betti numbers at each parameter and get the Betti curves $\beta_{0}^{a}$ and $\beta_{1}^{a}$. If the barcode obtained from the data set has a large number of bars, then using the persistent Betti number is often a reasonable choice.

				When it comes to persistent Laplacians, researchers seldom find it necessary to compute the persistent Laplacian $\Delta^{a,b}_{p}$ for all $a\leq b$ since it would make the calculation overly complex. In most cases, people only compute the Laplacians $\Delta^{a}_{0}$ and $\Delta^{a}_{1}$. It is yet to demonstrate that computationally sophisticated persistence offers any real benefit for data analysis.
			Therefore, we might misuse the term ``persistence'', which is largely due to computational considerations and does not affect the theoretical notion of persistence.
			 In our examples, we have not provided the ChatGPT code obtained with our supervision for computing the barcodes and $(a,b)$-persistent Laplacians for $a<b$. We only calculate the Betti numbers and Laplacians at various filtration parameters.

        \subsection{ChatGPT for persistent harmonic space }
	It is proved that the persistent harmonic space  for a filtration of simplicial complexes is naturally isomorphic to the persistent homology for the filtration \cite[Theorem 3.6]{liu2023algebraic}. It is worth noting that the basis for the persistent homology is a family of equivalent classes, which give  rise to an obstruction  for computing the generators in  the persistent homology. However, the persistent harmonic space  has a good property that the basis of the space is unique under linear isomorphisms, which provides us with a chance to find the basis of the persistent harmonic space and to track the generators as the persistence parameter varies.
				
 Although there are many existing  algorithms for computing persistent homology, there is not algorithm for computing the persistent harmonic space, to the best of our knowledge. It is interesting and important to demonstrate that ChatGPT can develop unseen codes in computational topology.

     Consider the simplicial complex $\mathcal{K}_{\varepsilon}$ parametrized by a real number $\varepsilon$ such that $\mathcal{K}_{\varepsilon}\hookrightarrow \mathcal{K}_{\varepsilon'}$ for any $\varepsilon\leq \epsilon'$. Recall that the $(a,b)$-persistent Laplacian $\Delta^{a,b}_{p}:C_{\ast}(\mathcal{K}_{a})\to C_{\ast}(\mathcal{K}_{a})$ is $\Delta^{a,b}_{p}=(d_{p}^{a})^{\ast}\circ d_{p}^{a}+(d_{p+1}^{a,b})^{\ast}\circ d_{p+1}^{a,b}$. The \emph{$(a,b)$-persistent harmonic space} is defined by
        \begin{equation*}
          \mathcal{H}^{a,b}_{p}=\ker \Delta^{a,b}_{p},\quad p\geq 0.
        \end{equation*}
        Note that there is a decomposition of $C_{\ast}(\mathcal{K}_{a})=\mathcal{H}^{a}_{\ast}\oplus \im d^{a}\oplus \im (d^{a})^{\ast}$, where $\mathcal{H}^{a}_{\ast}=\ker \Delta^{a}_{\ast}$. Let $p^{a}:C_{\ast}(\mathcal{K}_{a})\to \mathcal{H}^{a}_{\ast}$ be the projection of $C_{\ast}(\mathcal{K}_{a})$ to the harmonic part. By \cite[Proposition 3.8]{liu2023algebraic}, one has an isomorphism
        \begin{equation}\label{equation:harmonic}
          \mathcal{H}^{a,b}_{p}\cong p^{b}j_{\ast}^{a,b}(\mathcal{H}^{a}_{p})
        \end{equation}
        for any $a\leq b$. This implies the persistence on the persistent harmonic space. There is a natural isomorphism between the persistent harmonic space and the persistent homology
        \begin{equation*}
          \mathcal{H}^{a,b}_{p}\cong H^{a,b}_{p},\quad p\geq 0.
        \end{equation*}
        To compute the persistent homology, it suffices to give the code for the persistent harmonic space.

        Let $B_{p}^{a}$ be the representation matrix of the boundary operator $d_{p}^{a}:C_{p}(\mathcal{K}_{a})\to C_{p-1}(\mathcal{K}_{a})$. The representation matrix of the adjoint operator $(d_{p}^{a})^{\ast}:C_{p-1}(\mathcal{K}_{a})\to C_{p}(\mathcal{K}_{a})$ of $d_{p}^{a}$ is given by the transpose matrix $(B_{p}^{a})^{T}$ of $B_{p}^{a}$ with respect to the standard orthogonal basis given by the simplices. Thus, the harmonic space $\mathcal{H}_{p}^{a}$ can be represented by the solution space $W^{a}$ of the equation
        \begin{equation*}
          W^{a}\left(
                 \begin{array}{cc}
                   B_{p}^{a} & (B_{p+1}^{a})^{T} \\
                 \end{array}
               \right)=\mathbf{O}.
        \end{equation*}
        All the generators in $W^{a}$ have the index $a$ as the birth time. When the time comes to $b\geq a$, we will determine if there are generators that die at $b$. Note that $j_{\ast}^{a,b}$ is an inclusion, by Eq. (\ref{equation:harmonic}), we have that the $(a,b)$-persistent harmonic space $\mathcal{H}^{a,b}_{p}$ is the projection of $\mathcal{H}^{a}_{p}$ to $\mathcal{H}^{b}_{p}$. Note that $\mathcal{H}^{a}_{p}\subseteq \ker d^{b}_{p}$ and $\mathcal{H}^{b}_{p}=\ker d^{b}_{p}\cap \ker (d^{b}_{p+1})^{\ast}$. Let $U$ be the orthogonal complement space of $\mathcal{H}^{a,b}_{p}$ in $\mathcal{H}^{a}_{p}$. It follows that
        \begin{equation*}
         (d_{p+1}^{b})^{\ast}\mathcal{H}^{a}_{p}\cong U.
        \end{equation*}
        Thus, if $W^{a}(B_{p+1}^{b})^{T}=\mathbf{O}$, there is no generator dead at $b$. Moreover, $\dim W^{a}(B_{p+1}^{b})^{T}$ is the number of generators dead at $b$. The subspace of $W^{a}$
        corresponding to $W^{a}(B_{p+1}^{b})^{T}$ are the space generated by the dead generators. The left generators contribute to the persistent harmonic space $\mathcal{H}^{a,b}$, represented by the matrix $W^{a,b}$.

        Let $E^{a}_{p}$ be the basis of $C_{p}(\mathcal{K}_{a})$. Let $E^{b}_{p}=\left(
                     \begin{array}{c}
                       E^{a}_{p} \\
                       F^{b}_{p} \\
                     \end{array}
                   \right)$ be the basis of $C_{p}(\mathcal{K}_{b})$.
        Then we assume
        \begin{equation*}
          d_{p}^{b}\left(
                     \begin{array}{c}
                       E^{a}_{p} \\
                       F^{b}_{p} \\
                     \end{array}
                   \right)=\left(
                             \begin{array}{cc}
                               B_{p}^{a} & \mathbf{O} \\
                               \alpha_{p}^{b} & \beta_{p}^{b} \\
                             \end{array}
                           \right)\left(
                     \begin{array}{c}
                       E^{a}_{p-1} \\
                       F^{b}_{p-1} \\
                     \end{array}
                   \right).
        \end{equation*}
        Let $W_{b}=\left(
                     \begin{array}{cc}
                       W^{a,b} & \mathbf{O} \\
                       W^{b}_{1} & W^{b}_{2} \\
                     \end{array}
                   \right).
        $
        It follows that the the harmonic space at $b$ is given by $W^{b}\left(
                 \begin{array}{cc}
                   B_{p}^{b} & (B_{p+1}^{b})^{T} \\
                 \end{array}
               \right)=\mathbf{O}$. Thus one has
        \begin{equation}\label{equation:solution}
          \left(
                     \begin{array}{cc}
                       W^{b}_{1} & W^{b}_{2} \\
                     \end{array}
                   \right)\left(
                             \begin{array}{cccc}
                               B_{p}^{a} & \mathbf{O}&(B_{p+1}^{a})^{T}& (\alpha_{p+1}^{b})^{T}\\
                               \alpha_{p}^{b} & \beta_{p}^{b}&\mathbf{O}&(\beta_{p+1}^{b})^{T} \\
                             \end{array}
                           \right)=\mathbf{O}.
        \end{equation}
        Note that if $(\mathbf{x},\mathbf{0})$ is a solution vector, then $\mathbf{x}B_{p}^{a}=\mathbf{0}$ and $\mathbf{x}(\alpha_{p+1}^{b})^{T}=\mathbf{0}$. It follows that $\mathbf{x}$ is a generator representation in $\mathcal{H}_{p}^{a}$. Thus, all the generators represented in the form $(\mathbf{x},\mathbf{y})$ consist the solution space of Eq. (\ref{equation:solution}), where $\mathbf{y}\neq \mathbf{0}$. Therefore, one can obtain the desired generators and mark them with the birth time $b$. By continuing this iteration, we can obtain the persistent harmonic space. In particular, the corresponding dimensions for the persistent harmonic space lead to the persistence diagram or barcode for the persistent harmonic space, which is the consistent with the persistence diagram of the corresponding persistent homology.  In the example code of ChatGPT, we show the basic calculations for computing the 1-dimensional persistent harmonic space.

        \begin{figure}[H]
        \centering
        \begin{mdframed}[backgroundcolor=black, linewidth=0pt]
			\noindent\textcolor[rgb]{1.00,1.00,1.00}	{
				\texttt{\blue{import} vr\_complex\\
					\blue{import} numpy \blue{as} np\\
					\blue{import} boundary\_matrix\\
          			}}
		  \end{mdframed}
          \end{figure}
       \begin{figure}[H]
         \centering
         \begin{mdframed}[backgroundcolor=black, linewidth=0pt]
		 \noindent\textcolor[rgb]{1.00,1.00,1.00}	{
			\texttt{
					\blue{def} \red{generate\_complex}(A, B, d):\\
					\text{\qquad}result = [A[i] \blue{for} i, distance \blue{in} \orange{enumerate}(B) \blue{if} distance < d]\\
					\blue{return} result\\
                    \blue{def} \red{rips\_complex\_generate}(vr,a):\\
					\text{\qquad}rips\_complex\_1 = generate\_complex(vr[\red{1}],vr[\red{2}],a)\\
					\text{\qquad}rips\_complex\_2 = generate\_complex(vr[\red{3}],vr[\red{4}],a)\\
					\text{\qquad}rips\_complex = vr[\red{0}] + rips\_complex\_1 + rips\_complex\_2\\
					\blue{return} rips\_complex\\
					\\
					\blue{def} \red{dim\_index}(hypergraph):\\
					\text{\qquad}num\_hyperedges = \orange{len}(hypergraph)\\
					\text{\qquad}m=\orange{len}(hypergraph[num\_hyperedges-\red{1}])+\red{1}\\
					\text{\qquad}m=\orange{max}(m,\red{4})\\
					\text{\qquad}b = [\red{0}] * m\\
					\text{\qquad}\blue{for} i \blue{in} \orange{range}(num\_hyperedges):\\
					\text{\qquad}\text{\qquad}b[\orange{len}(hypergraph[i])]+=\red{1};  \\
					\text{\qquad}c = [\red{0}] * m\\
					\text{\qquad}\blue{for} i \blue{in} \orange{range}(m-\red{1}):\\
					\text{\qquad}\text{\qquad}c[i+\red{1}]=c[i]+b[i+\red{1}]\\
					\text{\qquad}\blue{return} c\\
					\\
					\blue{def} \red{left\_null\_space}(matrix):\\
					\text{\qquad}\gray{\# Calculate the rank of the input matrix}\\
					\text{\qquad}rank\_matrix = np.linalg.matrix\_rank(matrix)\\
					\\
					\text{\qquad}\gray{\# Perform SVD to get the left singular vectors}\\
					\text{\qquad}U, S, VT = np.linalg.svd(matrix)\\
					\\
					\text{\qquad}\gray{\# Extract the left null space vectors from U}\\
					\text{\qquad}left\_null\_space\_vectors = U[:, rank\_matrix:].T\\
					\\
					\text{\qquad}\blue{return} left\_null\_space\_vectors\\
					\\
					\blue{def} \red{triangular\_matrix}(Ar):\\
					\text{\qquad}\gray{\# Initialize the left multiplication matrix as the identity matrix}\\
					\text{\qquad}La = np.identity(Ar.shape[\red{0}])\\
					\text{\qquad}\gray{\# Iterate through each column and perform Gaussian elimination}
					\text{\qquad}\blue{for} j \blue{in} \orange{range}(Ar.shape[\red{1}]):\\
					\text{\qquad}\blue{for} i \blue{in} \orange{range}(j + \red{1}, Ar.shape[\red{0}]):\\
					\text{\qquad}\gray{\# Calculate the steps of Gaussian elimination to obtain elements of L}\\
					\text{\qquad}factor = Ar[i, j] / Ar[j, j]\\
					\text{\qquad}\gray{\# Update elements of L}\\
					\text{\qquad}La[i, j] = factor\\
					\text{\qquad}\gray{\# Update rows of A to achieve elimination}\\
					\text{\qquad}Ar[i, :] -= factor * Ar[j, :]\\
					\text{\qquad}\blue{return} La, Ar\\
         }}
		  \end{mdframed}
          \end{figure}
       \begin{figure}[H]
         \centering
         \begin{mdframed}[backgroundcolor=black, linewidth=0pt]
		 \noindent\textcolor[rgb]{1.00,1.00,1.00}	{
			\texttt{\gray{\# Calculate the inverse transformation of the matrix row operations for L}\\
					\blue{def} \red{lower\_triangular\_inverse}(A):\\
					\text{\qquad}n = A.shape[\red{0}]\\
					\text{\qquad}B = np.eye(n)\\
					\\
					\text{\qquad}\blue{if} np.linalg.matrix\_rank(A) < n:\\
					\text{\qquad}\text{\qquad}\orange{print}(\green{"The matrix is singular. Try another one!"})\\
					\text{\qquad}\text{\qquad}\blue{return} B\\
					\\
					\text{\qquad}\blue{for} k \blue{in} \orange{range}(n):\\
					\text{\qquad}\text{\qquad}B[\red{0}, k] = B[\red{0}, k] / A[\red{0}, \red{0}]\\
					\text{\qquad}\text{\qquad}\blue{for} i \blue{in} \orange{range}(\red{1}, n):\\
					\text{\qquad}\text{\qquad}\text{\qquad}B[i, k] = (B[i, k] - np.dot(A[i, :i], B[:i, k])) / A[i, i]\\
					\text{\qquad}\blue{return} B\\
					\\
					\gray{\#for a harmonic\_matrix, print the generators for the harmonic spaces}\\
					\blue{def} \red{print\_harmonic\_generators}(rips\_complex,harmonic\_matrix):\\
					\text{\qquad}\gray{\# Get the number of rows in the matrix}\\
					\text{\qquad}num\_rows = harmonic\_matrix.shape[\red{0}]\\
					\text{\qquad}num\_cols = harmonic\_matrix.shape[\red{1}]\\
					\text{\qquad}\gray{\#print the generators}\\
					\text{\qquad}\blue{for} i \blue{in} \orange{range}(num\_rows):\\
					\text{\qquad}\text{\qquad}\blue{for} j \blue{in} \orange{range}(num\_cols):\\
					\text{\qquad}\text{\qquad}\text{\qquad}\blue{if} j > \red{0}	 \blue{and} harmonic\_matrix[i][j]>\red{0} and j<num\_cols-\red{1}:\\
					\text{\qquad}\text{\qquad}\text{\qquad}\text{\qquad}\orange{print}(\green{"+"},round(harmonic\_matrix[i][j],\red{7}),rips\_complex[j],end=\green{""})\\
					\text{\qquad}\text{\qquad}\text{\qquad}\blue{elif} harmonic\_matrix[i][j] ==\red{0} and j<num\_cols-\red{1}:\\
					\text{\qquad}\text{\qquad}\text{\qquad}\text{\qquad}\blue{continue}\\
					\text{\qquad}\text{\qquad}\text{\qquad}\blue{elif} j<num\_cols-\red{1}:\\
					\text{\qquad}\text{\qquad}\text{\qquad}\text{\qquad}\orange{print}(round(harmonic\_matrix[i][j],\red{7}),rips\_complex[j],end=\green{""})\\
					\text{\qquad}\text{\qquad}\text{\qquad}\blue{elif} i<num\_rows-\red{1}:\\
					\text{\qquad}\text{\qquad}\text{\qquad}\text{\qquad}\orange{print}(\green{" and"})\\
					\text{\qquad}\text{\qquad}\text{\qquad}\blue{elif} harmonic\_matrix[i][j]>\red{0};\\
					\text{\qquad}\text{\qquad}\text{\qquad}\text{\qquad}\orange{print}(\green{"+"},round(harmonic\_matrix[i][j],\red{7}),rips\_complex[j])\\
					\text{\qquad}\text{\qquad}\text{\qquad}\blue{else}:\\
					\text{\qquad}\text{\qquad}\text{\qquad}\text{\qquad}\orange{print}(round(harmonic\_matrix[i][j],\red{7}),rips\_complex[j]))\\
					\text{\qquad}\blue{return}\\
					\\
					\gray{\# since not each matrix can be trangulated, we use the gaussian elimination}\\
					\blue{def} \red{gaussian\_elimination\_partial\_pivoting}(matrix):\\
					\text{\qquad}m, n = matrix.shape\\
					\text{\qquad}	\blue{for} i \blue{in} \orange{range}(\orange{min}(m, n)):\\
					\text{\qquad}\text{\qquad}\gray{\# Use partial pivoting to select the main element} \\
					\text{\qquad}\text{\qquad}max\_row = np.argmax(np.abs(matrix[i:, i])) + i\\
					\text{\qquad}\text{\qquad}\text{\qquad}\blue{if} matrix[max\_row, i] == \red{0}:\\
					\text{\qquad}\text{\qquad}\text{\qquad}\text{\qquad}\blue{continue}  \gray{\# If a non-zero main element cannot be found, move to the next column}	 \\
					\text{\qquad}\text{\qquad}\text{\qquad}matrix[[i, max\_row]] = matrix[[max\_row, i]]\\
					\text{\qquad}\text{\qquad}\text{\qquad}\blue{for} j \blue{in} \orange{range}(i + \red{1}, m):\\
					\text{\qquad}\text{\qquad}\text{\qquad}\text{\qquad}factor = matrix[j, i] / matrix[i, i]\\
					\text{\qquad}\text{\qquad}\text{\qquad}\text{\qquad}matrix[j, i:] -= factor * matrix[i, i:]\\
					\text{\qquad}\blue{return} matrix\\
				          			}}
		  \end{mdframed}
          \end{figure}
       \begin{figure}[H]
         \centering
         \begin{mdframed}[backgroundcolor=black, linewidth=0pt]
		 \noindent\textcolor[rgb]{1.00,1.00,1.00}	{
			\texttt{
				discrete\_points = np.array([[\red{0},\red{0},\red{0}],[\red{0},\red{1},\red{0}],\\
				\text{\qquad}\text{\qquad}\text{\qquad}\text{\qquad}\text{\qquad}\text{\qquad}\text{\qquad}[\red{1},\red{0},\red{0}],[\red{1},\red{1},\red{0}],\\
				\text{\qquad}\text{\qquad}\text{\qquad}\text{\qquad}\text{\qquad}\text{\qquad}\text{\qquad}[\red{0},\red{0},\red{2}],[\red{0},\red{1},\red{2}],\\
				\text{\qquad}\text{\qquad}\text{\qquad}\text{\qquad}\text{\qquad}\text{\qquad}\text{\qquad}[\red{1},\red{0},\red{2}],[\red{1},\red{1},\red{1}.\red{5}]\\
				\text{\qquad}\text{\qquad}\text{\qquad}\text{\qquad}\text{\qquad}\text{\qquad}\text{\qquad}])\\
				\\
				a = \red{1}.\red{1}\\
				b = \red{1}.\red{3}\\
				\gray{\#Give the maximal VR complex with all 0,1,2-dimensional simplices}\\
				vr\_list=vr\_complex.vr\_compute(discrete\_points)\\
				\\
				\gray{\#generate the rips complexes at some parameters}\\
				RipsC\_a = rips\_complex\_generate(vr\_list,a)\\
				RipsC\_b = rips\_complex\_generate(vr\_list,b)\\
				\\
				\orange{print}(\green{"The Rips complex at parameter a is:"},RipsC\_a)\\
				\orange{print}(\green{"The Rips complex at parameter b is:"},RipsC\_b)\\
                \\
				\gray{\#generate the boundary matrices for rips complexes}\\
				boundary\_a=boundary\_matrix.compute\_boundary\_matrix(RipsC\_a)\\
				boundary\_b=boundary\_matrix.compute\_boundary\_matrix(RipsC\_b)\\
				\\
				index\_a = dim\_index(RipsC\_a)\\
				index\_b = dim\_index(RipsC\_b)\\
				\\
				rips\_a = RipsC\_a[index\_a[\red{1}]:index\_a[\red{2}]]\\
				rips\_b = RipsC\_b[index\_b[\red{1}]:index\_b[\red{2}]]\\
				\\
				\\
				\gray{\#compute W\_a, which is the harmonic matrix at a, }\\
				\gray{\#where dirac\_1 is the union of boundary and coboundy matrices}\\
				dirac\_1 = np.hstack((boundary\_a[\red{0}], boundary\_a[\red{1}].T))\\
				W\_a = left\_null\_space(dirac\_1)\\
				\\
				\gray{\#if W\_a is emtpyset, the persistent harmonic space is none, or we compute the matrix U in our paper}\\
				\blue{if} W\_a.size == \red{0}:\\
				\text{\qquad}\orange{print}(\green{"The harmonic space at a is none."})\\
				\text{\qquad}\orange{print}(\green{"Then the harmonic space from a to b is also none."})\\
				\blue{else}:\\
				\text{\qquad}\orange{print}(\green{"The harmonic space at a is generated by:"})\\
				\text{\qquad}print\_harmonic\_generators(rips\_a,W\_a)\\
				\\
				\text{\qquad}\gray{\#U\_matrix is the space of dead elements at b}\\
				\\
				\text{\qquad}\gray{\#The number of one-dimensional simplices added in Rips complex at time b compared to time a.}\\
				\text{\qquad}rank\_2 = index\_b[\red{2}]-index\_b[\red{1}]-index\_a[\red{2}]+index\_a[\red{1}]\\
                           			}}
		  \end{mdframed}
          \end{figure}
       \begin{figure}[H]
         \centering
         \begin{mdframed}[backgroundcolor=black, linewidth=0pt]
		 \noindent\textcolor[rgb]{1.00,1.00,1.00}	{
			\texttt{
			    \text{\qquad}\gray{\#The number of points}\\
				\text{\qquad}rank\_1 = W\_a.shape[\red{0}]\\
				\text{\qquad}zero\_add = np.zeros((rank\_1, rank\_2))\\
				\\
				\text{\qquad}W\_a\_add = np.hstack((W\_a,zero\_add))\\
				\\
				\text{\qquad}U\_matrix = W\_a\_add@(boundary\_b[\red{1}].T)\\
				\text{\qquad}\blue{if} U\_matrix.size == \red{0}:\\
				\text{\qquad}\text{\qquad}\orange{print}(\green{"The harmonic space from a to b is generated by:"})\\
				\text{\qquad}\text{\qquad}W\_ab = W\_a\_add\\
				\text{\qquad}\text{\qquad}print\_harmonic\_generators(rips\_b,W\_ab)\\
				\text{\qquad}\blue{else}:\\
				\text{\qquad}\text{\qquad}rank\_U = np.linalg.matrix\_rank(U\_matrix)\\
				\text{\qquad}\text{\qquad}\blue{if}  rank\_U == rank\_1:\\
				\text{\qquad}\text{\qquad}\text{\qquad}\orange{print}(\green{"The harmonic space from a to b is none."})\\
				\text{\qquad}\text{\qquad}\blue{else}: \\
				\text{\qquad}\text{\qquad}\text{\qquad}\gray{\# compute the persistent harmonic space matrix W\_ab from a to b}\\
				\text{\qquad}\text{\qquad}\text{\qquad}Left\_triangle,trangle\_matrix = triangular\_matrix(U\_matrix)\\
				\text{\qquad}\text{\qquad}\text{\qquad}copy\_triangle= Left\_triangle[:].copy()\\
				\text{\qquad}\text{\qquad}\text{\qquad}L\_inv = lower\_triangular\_inverse(copy\_triangle)\\
				\text{\qquad}\text{\qquad}\text{\qquad}W\_temp = L\_inv @ W\_a\\
				\text{\qquad}\text{\qquad}\text{\qquad}W\_ab = W\_temp[:rank\_U, :]\\
				\\
				\gray{\#W\_b is the harmonic matrix at b}\\
				dirac\_2 = np.hstack((boundary\_b[\red{0}], boundary\_b[\red{1}].T))\\
				W\_b = left\_null\_space(dirac\_2)\\
                \\
				\blue{if} W\_b.size == 0:\\
				\text{\qquad}\orange{print}(\green{"The harmonic space at b is none."})\\
				\blue{else}:\\
				\text{\qquad}\orange{print}(\green{"The harmonic space at b is generated by:"})\\
				\text{\qquad}print\_harmonic\_generators(rips\_b,W\_b)\\
				\text{\qquad}\_, birth\_b = triangular\_matrix(W\_b)\\
				\text{\qquad}reflect\_matrix = np.fliplr(birth\_b)\\
				\text{\qquad}upper\_triangular\_matrix = gaussian\_elimination\_partial\_pivoting (reflect\_matrix)\\
				\text{\qquad}rr\_matrix = np.fliplr(upper\_triangular\_matrix)\\
				\text{\qquad}\blue{if} W\_a.size == 0:\\
				\text{\qquad}\text{\qquad}\orange{print}(\green{"The harmonic generators birth at b is generated by:"})\\
				\text{\qquad}\text{\qquad}print\_harmonic\_generators(rips\_b,W\_b)\\
				\text{\qquad}\blue{else}	:\\
				\text{\qquad}\text{\qquad}number\_of\_generators = W\_b.shape[\red{0}] - W\_ab.shape[\red{0}]\\
				\text{\qquad}\text{\qquad}\orange{print}(\green{"The harmonic generators born at b is generated by:"})\\
				\text{\qquad}\text{\qquad}print\_harmonic\_generators(rips\_b,W\_b[:number\_of\_generators, :]) \\
			}}
		\end{mdframed}
        \caption{There are ten functions appeared in the ChatGPT code. The functions \textbf{compute\_boundary\_matrix} and \textbf{vr\_compute} are imported from the previous files \textbf{boundary\_matrix} and \textbf{vr\_complex}, respectively. Apart from the function \textbf{print\_harmonic\_generators}, the functions \textbf{rips\_complex\_generate} and \textbf{dim\_index} were designed by ChatGPT with our guidance. All the other functions, including \textbf{generate\_complex}, \textbf{left\_null\_space}, \textbf{triangular\_matrix}, \textbf{lower\_triangular\_inverse}, and \textbf{gaussian\_elimination\_partial\_pivoting}, were provided in responses by ChatGPT. Specifically, the main idea for the function \textbf{gaussian\_elimination\_partial\_pivoting} was also provided by ChatGPT.}\label{code_harmonic_space}
        \end{figure}

        \begin{example}
        Consider the point set $X$ embedded into the Euclidean space $\mathbb{R}^{3}$ as follows.
        \begin{equation*}
          X=\{(0,0,0),(0,1,0),(1,0,0),(1,1,0),(0,0,2),(0,1,2),(1,0,2),(1,1,1.5)\}.
        \end{equation*}
        We guided ChatGPT to complete the code for computing the 1-dimensional persistent harmonic space. We take the persistent parameter $a=1.2$ and $b=1.3$. Then the result window displays as follows.
        \begin{figure}[H]
        \begin{small}
			\centering
			\begin{mdframed}[backgroundcolor=blue!30, linewidth=0pt]
				\noindent   The Rips complex at parameter a is: [[0], [1], [2], [3], [4], [5], [6], [7], [0, 1], [0, 2], [1, 3], [2, 3], [4, 5], [4, 6]]\\
          The Rips complex at parameter b is:
          [[0], [1], [2], [3], [4], [5], [6], [7], [0, 1], [0, 2], [1, 3], [2, 3], [4, 5], [4, 6], [5, 7], [6, 7]]\\
          The harmonic space at a is generated by:
          0.5 [0, 1]-0.5 [0, 2]+ 0.5 [1, 3]-0.5 [2, 3]\\
          The harmonic space from a to b is generated by:
          0.5 [0, 1]-0.5 [0, 2]+ 0.5 [1, 3]-0.5 [2, 3]\\
          The harmonic space at b is generated by:
          -0.5 [0, 1]+ 0.5 [0, 2]-0.5 [1, 3]+ 0.5 [2, 3] and
          -0.5 [4, 5]+ 0.5 [4, 6]-0.5 [5, 7]+ 0.5 [6,7]\\
          The harmonic generators born at b is generated by:
          -0.5 [4, 5]+ 0.5 [4, 6]-0.5 [5, 7]+0.5 [6,7]\\
			\end{mdframed}
          \end{small}
			\caption{Shown results on persistent harmonic spaces (I).}\label{fig:harmonic_1}
		\end{figure}
        Besides, for another example, $X$ is the point set with the points $(0.5, \sqrt{3}/ 2, 0)$, $(0.5, \sqrt{3} / 2, 0)$, $(-1, 0, 0)$ ,$(1, 0, 0)$, $(0.5, -\sqrt{3} / 2, 0)$ ,$(-0.5, -\sqrt{3} / 2, 0.3)$. We take $a=1$ and $b=1.2$. The result is shown in Fig. \ref{fig:harmonic_2}.
        \begin{figure}[H]
        \begin{small}
			\centering
			\begin{mdframed}[backgroundcolor=blue!30, linewidth=0pt]
				\noindent   The Rips complex at parameter a is: [[0], [1], [2], [3], [4], [5], [0, 3], [1, 2], [3, 4]]\\
          The Rips complex at parameter b is: [[0], [1], [2], [3], [4], [5], [0, 3], [1, 2], [3, 4], [0, 1], [2, 5], [4, 5]]\\
          The harmonic space at a is none.\\
          Then the harmonic space from a to b is also none.\\
          The harmonic space at b is generated by:\\
          -0.4082483 [0, 3]+ 0.4082483 [1, 2]-0.4082483 [3, 4]+ 0.4082483 [0, 1]+ 0.4082483 [2, 5]\\
          The harmonic generators birth at b is generated by:\\
          -0.4082483 [0, 3]+ 0.4082483 [1, 2]-0.4082483 [3, 4]+ 0.4082483 [0, 1]+ 0.4082483 [2, 5]
			\end{mdframed}
          \end{small}
			\caption{Shown results on the persistent harmonic space (II).}\label{fig:harmonic_2}
		\end{figure}
		
        In the example, one see that we can indeed obtain the persistent harmonic space, the corresponding harmonic generators, and the generators of birth and death.
        If we only consider the dimensions of harmonic spaces, the code can be developed into a code for computing the barcode for persistent dimensions of harmonic spaces.
        \end{example}
				
        Through our interactions with ChatGPT, we have discussed various aspects of matrix transformations, operations, and calculations. We have gained knowledge about matrix coding, and ChatGPT has improved its understanding of the questions we have posed during our Q\&A sessions.

        As demonstrated above, we can guide ChatGPT step by step to implement the required code. What we want to emphasize is that ChatGPT can not only help us find existing codes but also  create new codes that are unknown to researchers.       This highlights the creative potential of ChatGPT for computations.

		\section{Applications}\label{section:application}

		In this section, the ChatGPT codes obtained with our supervision are used to calculate data sets as examples. We will present two examples: one involving a point set of six points arranged as nodes of a hexagon, and the other involving   three molecules.
		
		\begin{figure}[H]
			\centering
			\includegraphics[width=12cm]{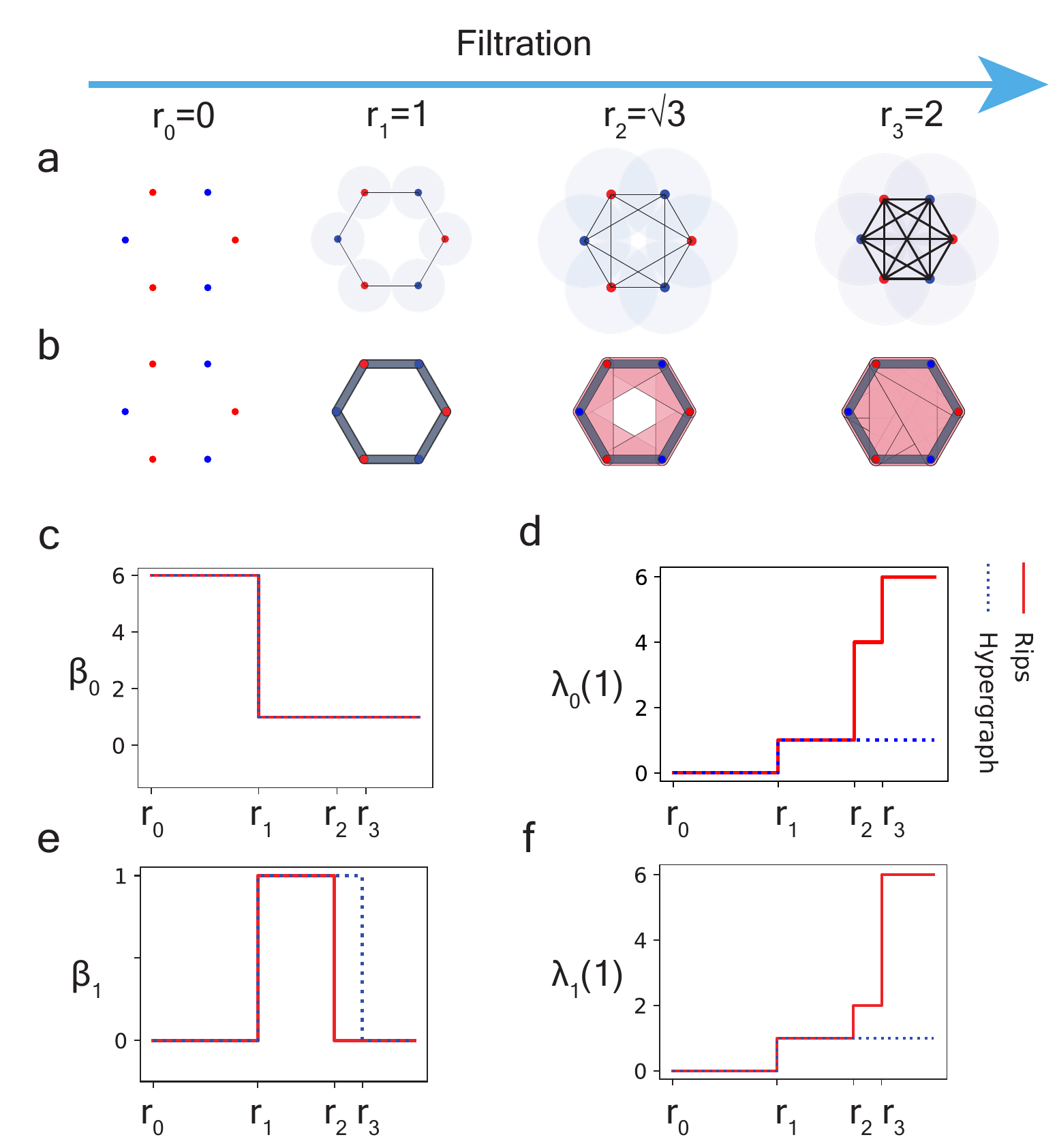}\\
			\caption{Comparison of persistent Laplacians generated by Rips complexes and hypergraphs.
			{\bf a} The filtration of Rips complex. {\bf b} The filtration of hypergraphs. The source code for generating virtualization of hypergraphs is from \cite{chen2023persistent}. {\bf c} Comparison of 0-dimensional Betti numbers for filtrations of Rips complexes and hypergraphs. {\bf d} Comparison of spectral gaps of 0-dimensional Laplacians for filtrations of Rips complexes and hypergraphs. {\bf e} Comparison of 1-dimensional Betti numbers for filtrations of Rips complexes and hypergraphs. {\bf f} Comparison of spectral gaps of 1-dimensional Laplacians for filtrations of Rips complexes and hypergraphs.}\label{example:hexagon}
		\end{figure}
		In Fig. \ref{example:hexagon}, we consider a point set composed of six vertices forming a regular hexagon with side lengths equal to 2. The six points are alternately colored in red and blue. For simplicial complex model, we use the Rips complex. For hypergraph model, the hyperedges are constructed from the corresponding Rips complex by requiring that each hyperedge has two colors. Specifically, the vertices are considered as 0-dimensional hyperedges. When the filtration parameter reaches $r_{1}=1$, the Rips complex forms a regular hexagon. The corresponding hypergraphs is given by $\mathcal{H}=\{\{0\},\{1\},\{2\},\{3\},\{4\},\{5\},\{0,1\},\{1,2\},\{2,3\},\{3,4\},\{4,5\},\{0,5\}\}$. Here, $0,\dots,5$ denotes the six points for convenience. When the filtration parameter $r_{2}=\sqrt{3}$, the points $0$ and $2$ are connected by an edge. But $\{0,2\}$ is not a hyperedge since $0$ and $2$ have the same color. Similarly, $\{1,3\}$, $\{2,4\}$, $\{3,5\}$, $\{0,4\}$, $\{1,5\}$ are not hyperedges. However, sets as $\{0,1,2\}$, $\{1,2,3\}$, $\{2,3,4\}$, $\{3,4,5\}$, $\{0,4,5\}$, $\{0,1,5\}$ are 2-dimensional hyperedges. When $r_{3}=2$, all the six points are connected. The corresponding hypergraph have more 2-dimensional hyperedges, such as $\{0,1,3\},\{0,1,4\}$ and so on. In \ref{example:hexagon}{\bf c}, the 0-th Betti numbers for the Rips complex model and the hypergraph model coincides. Fig. \ref{example:hexagon}{\bf e} shows difference of the 1-th Betti numbers between the two models. Fig. \ref{example:hexagon}{\bf d} and Fig. \ref{example:hexagon}{\bf f} show the spectral gaps of 0 and 1-dimensional Laplacians for the Rips complex model and the hypergraph model. When the spectral gap equals 0, it means that all the eigenvalues are zero. While the Rips complex model and the hypergraph give different information for the given   set of point cloud, the Rips complex model apparently offers rich information than a simple  hypergraph does. However, hypergraph can be tailored to reveal specific structural features.

			\begin{figure}[H]
			\centering
			\includegraphics[width=16cm]{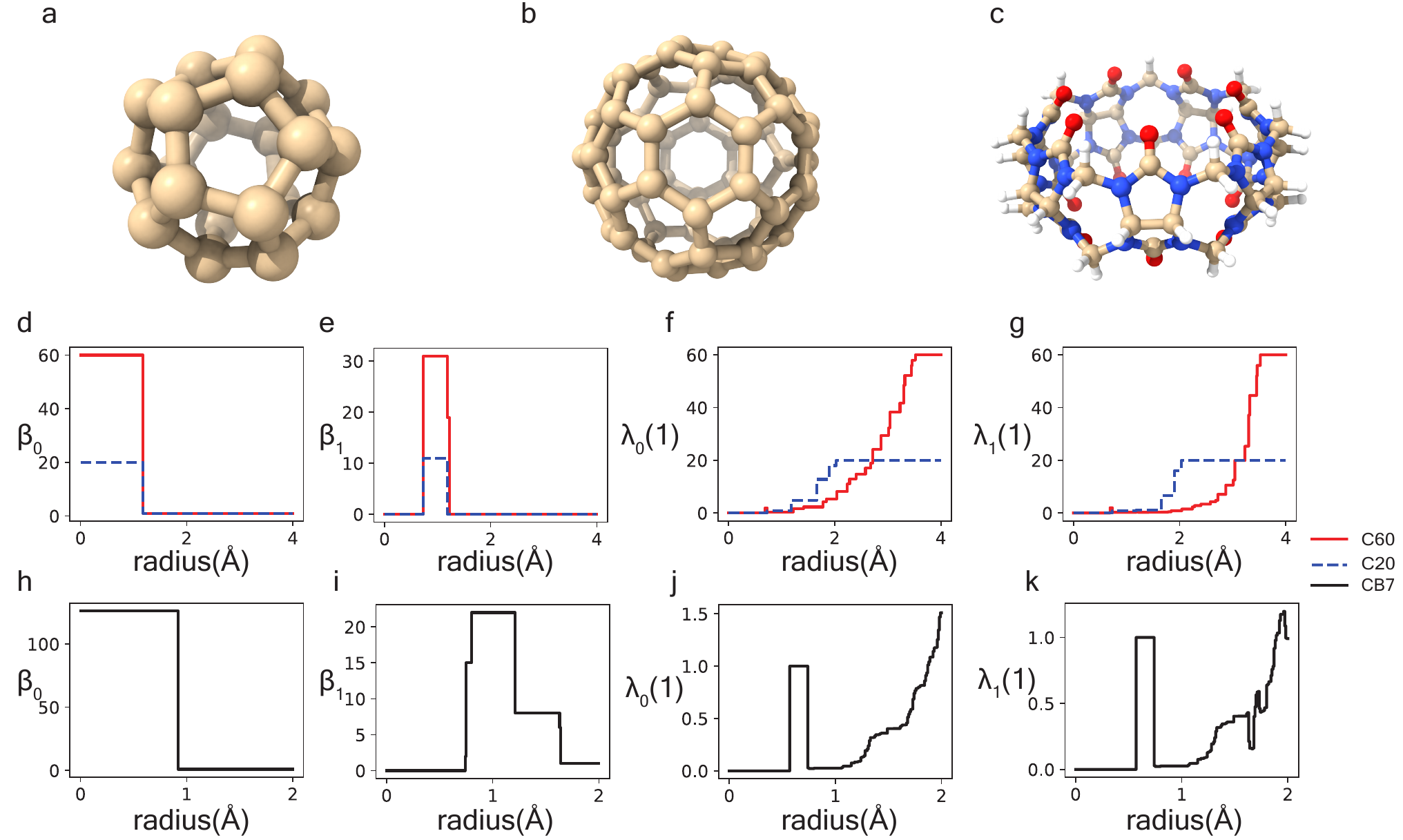}\\
			\caption{Illustration of persistent homology and persistent Laplacians for three molecules.
			{\bf a} The structure of fullerene $\mathrm{C}_{20}$. {\bf b} The structure of fullerene $\mathrm{C}_{60}$. {\bf c} The structure of cucurbit[7]uril $\mathrm{CB}7$. {\bf d} Comparison of 0-dimensional Betti numbers for $\mathrm{C}_{20}$ and $\mathrm{C}_{60}$. {\bf e} Comparison of 1-dimensional Betti numbers for $\mathrm{C}_{20}$ and $\mathrm{C}_{60}$. {\bf f} Comparison of spectral gaps of 0-dimensional Laplacians for $\mathrm{C}_{20}$ and $\mathrm{C}_{60}$. {\bf g} Comparison of spectral gaps of 1-dimensional Laplacians for $\mathrm{C}_{20}$ and $\mathrm{C}_{60}$. {\bf h} The 0-dimensional Betti numbers for $\mathrm{CB}7$. {\bf i} The 1-dimensional Betti numbers for $\mathrm{CB}7$. {\bf j} The spectral gaps of 0-dimensional Laplacians for $\mathrm{CB}7$. {\bf k} The spectral gaps of 1-dimensional Laplacians for $\mathrm{CB}7$.
			}\label{example:biomolecules}
		\end{figure}

		Additionally, we   consider examples of three  molecules, including $\mathrm{C}_{20}$, $\mathrm{C}_{60}$, and $\mathrm{BC}7$, and calculate the corresponding Betti numbers and Laplacian spectra over filtration.
	Note that $\mathrm{C}_{20}$ and $\mathrm{C}_{60}$ were 	studied in the literature using persistent   Laplacians \cite{wang2020persistent}, and  $\mathrm{BC}7$  was examined  in the literature using persistent path Laplacians \cite{wang2023persistent}.
		The topological features are computed from the Rips complexes of these small molecules. The atoms are regarded as points in Euclidean space and filtration lengths are determined by the distance between points.
	
		In Fig. \ref{example:biomolecules},  subfigures {\bf d}, {\bf e}, {\bf f}, and {\bf g} show the comparison of Betti numbers and the smallest positive eigenvalues for the corresponding Laplacians of dimension 0,1 for $\mathrm{C}_{20}$ and $\mathrm{C}_{60}$.
		Clearly, information obtained from persistent homology is limited. In contrast, persistent Laplacians not only return all the topological invariants as given by persistent homology in their harmonic spectra, but also reveal additional structural changes during the filtration in their nonharmonic spectra.
	Subfigures {\bf h}, {\bf i}, {\bf j}, and {\bf k} give the Betti numbers and the smallest positive eigenvalues for the corresponding Laplacians of dimension 0 and 1 of $\mathrm{BC}7$. The curves of Betti numbers and the smallest positive eigenvalues based on persistent Laplacians are   quite similar with those obtained with persistent path Laplacians \cite{wang2023persistent}.

\section{Concluding remarks}

In this study, we delve into the application of ChatGPT in computational topology from the perspective of a topologist who possesses limited knowledge of computational algorithms and lacks coding experience. To enable ChatGPT to generate accurate codes for computational topology tasks, we took several strategic steps. These steps involved training ChatGPT to understand fundamental mathematical concepts, guiding its responses in the right direction, leveraging ChatGPT to enhance our comprehension of the resulting algorithms, and designing test examples to validate the code provided by ChatGPT. This approach essentially involves a form of expert-supervised ChatGPT debugging.

One of the central mathematical models we focused on is the Vietoris-Rips complex, a widely used tool in applied topology that excels at characterizing the high-dimensional structure of data. Initially, we concentrated on algorithms for computing Betti numbers, Laplacians, and Dirac operators for simplicial complexes. Notably, the recently developed persistent Laplacians have greatly improved the effectiveness of applying topology in the field of molecular biology.

While simplicial complex models suit most applications, there are scenarios where they may not fully capture the data due to their  constraints. This is where hypergraph and directed graph models come into play. Both hypergraphs and directed graphs exhibit high-dimensional structures, with hypergraphs characterized by embedded homology and directed graphs characterized by path homology. In our study, we harnessed the capabilities of ChatGPT to implement codes for computing Betti numbers and Laplacians specifically tailored for hypergraphs and directed graphs. We also demonstrate ChatGPT's ability to compute the persistent harmonics space, which has been computed to the best of our knowledge.

By combining the expertise of a topologist with the capabilities of ChatGPT, we aim to advance the field of computational topology and make these powerful mathematical tools more accessible and applicable in various domains.

In our applications, we examined a point set comprising six vertices arranged in a regular hexagon. We carried out computations for the Vietoris-Rips complexes and determined the 0th and 1st-dimensional Betti numbers, as well as the smallest positive eigenvalue of Laplacians for hypergraphs. Subsequently, we generated plots illustrating the Betti curves and curves depicting the smallest positive eigenvalues.
The ChatGPT-generated code yielded accurate results in these calculations. Additionally, we extended our analysis to include three  molecules $\mathrm{C}{20}$, $\mathrm{C}{60}$, and $\mathrm{CB}7$. For these molecules, we computed their 0-dimensional and 1-dimensional Betti numbers, along with the smallest positive eigenvalue of their Laplacians. We then visualized the corresponding Betti curves and smallest positive eigenvalue curves. The efficiency of the code was reasonable, with the calculations taking approximately several tens of minutes to run on a standard laptop.

		In the future, there remain several unexplored avenues for further work and research.
Firstly, one critical aspect we have not addressed in this work is algorithm optimization. On one hand, ChatGPT could directly aid us in optimizing algorithms if desired. On the other hand, we could leverage ChatGPT to explore mathematical optimizations of algorithms, ultimately achieving more efficient computations.

Secondly, while the algorithms presented in this paper can calculate persistent Betti numbers and persistent Laplacians, they do not encompass the computation of barcodes. We believe that ChatGPT can play a pivotal role in helping us complete this task.

Lastly, it is worth noting that commonly used topological invariants include (co)homology, Laplacian, and Dirac. However, there are other algorithms for topological invariants, such as cup-product, Steenrod algebra, K-theory, and others, which warrant further investigation and study.
		
			This work signifies merely the initial stride in the utilization of chatbots for computational topology. We are confident that this step marks the commencement of a new era in mathematical research, forging connections between pure mathematical theories and computational tools in applied sciences. We envisage a promising future for this approach, as generative AI becomes increasingly potent and its capacity to comprehend mathematical concepts continuously enhances rapidly.

		\section*{Data and Code Availability}
		The data and source code obtained in this work are publicly available in the Github repository: \url{https://github.com/JoybearLiu/ChatGPT-for-computational-topology}.

		\section*{Acknowledgments}
		This work was supported in part by NIH grants R01GM126189, R01AI164266, and R35GM148196, National Science Foundation grants DMS2052983, DMS-1761320, and IIS-1900473, NASA  grant 80NSSC21M0023,   Michigan State University Research Foundation, and  Bristol-Myers Squibb  65109.

		\bibliographystyle{alpha}
		\bibliography{Reference}

	\end{CJK*}
\end{document}